\documentclass[a4,12pt]{amsart} 

\usepackage[utf8]{inputenc}
\usepackage[T1]{fontenc}
\usepackage{yfonts}
\usepackage[english]{babel}
\usepackage[normalem]{ulem}
\usepackage{cancel}
\usepackage{bm}

\usepackage{lipsum}
\usepackage{amsmath}
\usepackage{amsthm}
\usepackage{amssymb}

\usepackage[shortlabels]{enumitem}
\usepackage{graphicx}
\usepackage{mathtools}
\usepackage{hyperref}
\usepackage{amsfonts}
\usepackage{latexsym}
\usepackage{amscd}

\usepackage[all]{xy}
\usepackage[dvipsnames]{xcolor}

\def\cdeg{\deg_c}

\definecolor{light-gray}{gray}{0.55}

\newtheorem{theorem}{Theorem}[section]
\newtheorem{lemma}[theorem]{Lemma}
\newtheorem{proposition}[theorem]{Proposition}

\newtheorem{definition}[theorem]{Definition}
\theoremstyle{definition}
\newtheorem{remark}[theorem]{Remark}
\newtheorem{example}[theorem]{Example}
\begin{document}

\title[Injectives  over Leavitt path algebras]{Injectives  over Leavitt path algebras of graphs with disjoint cycles}
\author{Gene Abrams
}
\address{Department of Mathematics, University of Colorado,
Colorado Springs, CO 80918 U.S.A.}
\email{abrams@math.uccs.edu}
\thanks{
The second and third authors are supported by Progetto PRIN ``Categories, Algebras: Ring-Theoretical and Homological Approaches (CARTHA) '', chief advisor Alberto Facchini (2019 - 2022).   Part of this work was carried out during a visit of the third author to the University of Colorado  Colorado Springs.    The third author is pleased to take this opportunity to again express his thanks to the host institution, and its faculty, for its warm hospitality and support.  \\
}
\author{Francesca Mantese}
\address{Dipartimento di Informatica, Universit\`{a} degli Studi di Verona, I-37134 Verona, Italy}
\email{francesca.mantese@univr.it}

\author{Alberto Tonolo}
\address{Dipartimento di Matematica ``Tullio Levi-Civita'', Universit\`{a} degli Studi di Padova, I-35121, Padova, Italy}
\email{alberto.tonolo@unipd.it}

%
%
%
%
%

%


%
\maketitle

%

\begin{abstract}
Let $K$ be any field, and let $E$ be a finite graph with the property that every vertex in $E$ is the base of at most one cycle ({i.e., a graph with disjoint cycles}).  We explicitly construct the injective envelope of each simple left module over the Leavitt path algebra  $L_K(E)$.  
The main idea girding our construction is that of a ``formal power series'' extension of modules, thereby  developing for all graphs {with disjoint cycles} the understanding of injective envelopes of simple modules over $L_K(E)$ achieved previously for the simple modules over the Toeplitz algebra.


\normalsize

\end{abstract}


%
%

\section{Introduction}

In \cite{AMT21}, the three authors presented a complete, explicit description of the injective envelope of every simple left module over  $L_K(\mathcal{T})$, the ``Toeplitz Leavitt path $K$-algebra'' associated to the ``Toeplitz graph''  $\mathcal{T}:= \hspace{-.25in} \xymatrix{   & {\bullet} \ar@(ul,dl) \ar[r] &{\bullet}}$, {isomorphic to the well-studied Jacobson algebra $K\langle X,Y  \ | \ XY=1 \rangle$  (see for instance \cite{B, G, I, J}).}
Clearly $\mathcal{T}$ is an elementary but nontrivial example of a graph that  satisfies the following property.  

\begin{definition}\label{ARdefinition}
Let $E$ be a finite graph.  We say that $E$  {\rm {has disjoint cycles}} in case each vertex of $E$ is the base of at most one cycle.
\end{definition}


The condition that each vertex in a finite graph $E$ is the base of at most one cycle  has been studied in other settings.  For example, in  \cite{AAJZ} the four authors establish that, for any field $K$,  this condition is equivalent to the Leavitt path algebra $L_K(E)$ having finite Gelfand-Kirillov dimension.  Subsequently, Ara and Rangaswamy in \cite[Theorem 1.1]{AR14} establish that this condition is also equivalent to the property that all simple modules over $L_K(E)$ are of a specified type (so-called {\it Chen} modules).
{More recently, in \cite{HSV22} these graphs have been characterised in terms of the Jordan-H\"older composition series of their talented monoid.}

For us, {the condition to have disjoint cycles} sits at a wonderfully fortunate confluence of the ideas presented in \cite{AAJZ} and \cite{AR14}.   On  one hand, 
similar to a process carried out  in \cite{AAJZ}, 
 the condition affords the possibility to consider a well-defined positive integer representing  the  ``length of a chain of cycles'', and in particular affords the possibility of applying an induction argument in such a situation.   On the other hand, as provided by \cite{AR14},  having a complete explicit description of all the simple $L_K(E)$-modules puts us in position to describe the injective envelopes of all simple modules, and consequently to describe an injective cogenerator for $L_K(E)$-Mod.  

 We show in the current article that the ``formal power series'' idea developed in \cite[Section 6]{AMT21} in the context of the graph $\mathcal{T}$ can indeed be  extended to all graphs {with disjoint cycles}.   In particular, we show that the injective envelopes of all the simple modules over $L_K(E)$ may be realized  as a type of  formal power series extension.   This is the heart of our main result, Theorem \ref{thm:main}. 
 {The proof,  rather long and  delicate, is obtained first reducing via Morita equivalences to connected graphs which contain no \emph{source vertices} and in which every \emph{source cycle} is a loop, and then on an induction reasoning on the number of cycles.}  {The presence of cycles in the graph gives rise to  elements of $L_K(E)$ which can be seen as polynomials with the cycles as indeterminate. It is therefore not surprising to come across formal series when dealing with this class of algebras.}
  
The article is organized as follows.   In Section \ref{foundation} we review some of the basics of Leavitt path algebras, and provide some properties of Leavitt path algebras of graphs having a specified structure.  
 In Section \ref{section:ChenPrufer} we review the notion of a Chen simple $L_K(E)$-module,  
  as well as that of a Pr\"{u}fer $L_K(E)$-module.   
 In Section \ref{fps} we introduce the formal power series construction associated to subsets of ${\rm Path}(E)$.
With the Sections \ref{section:ChenPrufer} and \ref{fps} material in hand, in Section \ref{mainsection} we state  and prove the main result of the article, Theorem \ref{thm:main}. 

We conclude this introductory section by giving a reformulation of the standard Baer Criterion for injectivity, one which is well suited to our situation.

\begin{proposition}\label{injectivebystep}
Let $L$ be an associative ring, and let $M$ be a left $L$-module.  Let $I$ be a fixed left ideal of $L$. 
Assume that any morphism $f: X \to M$ from any left ideal $X \leq I$ extends to a morphism $\hat f:I\to M$, and also assume that any morphism $g: Y\to M$ from any left ideal $Y \geq I$ extends to a morphism $\hat g:L\to M$.
Then $M$ is injective.
\end{proposition}
\begin{proof}
Let $J$ be any left ideal of $L$, and let $h:J\to M$ be an $L$-module morphism. The restriction $h_0:J\cap I\to M$ extends by assumption to $\hat h_0:I\to M$. Setting  
$$h_1(j+i)=h(j)+\hat h_0(i)$$
for each $j\in J, i\in I$  yields  a morphism $h_1:J+I\to M$. That $h_1$ is well-defined follows from this observation:  $j+i=j'+i'$ implies $j-j'=i'-i\in J\cap I$, and hence 
\[h(j)-h(j')=h(j-j')=h_0(j-j')=h_0(i'-i)=\hat h_0(i'-i)=\hat h_0(i')-\hat h_0(i),\]
from which $h(j)+\hat h_0(i)=h(j')+\hat h_0(i')$. Clearly $h_1$ extends $h$.  Now $h_1:J+I\to M$ extends by assumption to $\hat h_1:L\to M$. By the Baer Criterion we conclude that $M$ is injective.
\end{proof}

\section{Leavitt path algebras}\label{foundation}

{A (directed)  graph $E=(E^0,E^1,r,s)$ consists of two  sets $E^0$ and $E^1$ together with maps $r,s:E^1\to E^0$, the \emph{range} and the \emph{source} maps. The elements of $E^0$ are called \emph{vertices}, those of $E^1$ \emph{edges}.   $E$ is \emph{finite} in case both $E^0$ and $E^1$ are finite sets.  {\bf Unless otherwise indicated, we assume throughout that all graphs are finite.}  A \emph{sink} $w$ is a vertex which emits no edges, {i.e., $s^{-1}(w)=\emptyset$; a \emph{source} $u$ is a vertex  to which no edges arrive, i.e., $r^{-1}(u)=\emptyset$}. A \emph{finite path} of length $n$ is a sequence of edges $\rho=e_1e_2\cdots e_n$ with $r(e_i)=s(e_{i+1})$ for $i=1,2,...,n-1$.  We denote $s(\rho) = s(e_1)$ and $r(\rho) = r(e_n)$. Any vertex is viewed as a (trivial) path of length 0.     We denote by $\text{Path}(E)$ the set of all finite paths in $E$. For each $e\in E^1$, we call $e^*$ the associated \emph{ghost edge}; by definition,  $r(e^*)=s(e)$ and $s(e^*)=r(e)$.

Given any field $K$ and (finite) graph $E$, the \emph{Leavitt path $K$-algebra} $L_K(E)$ is the free associative $K$-algebra generated by a set of symbols  $\{v:v\in E^0\} \sqcup \{e,e^*:e\in E^1\}$, where  $\{v:v\in E^0\}$ are orthogonal idempotents, and for which the following relations are imposed:
\begin{enumerate}
\item for each $e\in E^1$, \ $s(e)e=e=er(e)$ and $r(e)e^*=e^*=e^*s(e)$;
\item  for each pair $e$, $f\in E^1$, \ $e^*f=\begin{cases}
    r(f)  & \text{if }e=f, \\
      0& \text{otherwise} ;   \ \ \ \mbox{and}
\end{cases}$
\item for each non-sink $v\in E^0$, \ $v=\sum_{e\in s^{-1}(v)}ee^*$.
\end{enumerate}
The elements of $L_K(E)$ are $K$-linear combinations of paths $\lambda \mu^*$, $\lambda,\mu\in\text{Path}(E)$ \cite[Lemma~1.2.12]{AAS17}.
A nontrivial path $e_1e_2\cdots e_n$ is
a \emph{closed path} (resp., a \emph{cycle}\footnote{however, see Definition \ref{cycledef} below}) if $r(e_n)=s(e_1) $ (resp., and $s(e_i)\not=s(e_j)$ for every $i\not=j$). A cycle of length $1$ is called a \emph{loop}.  
A \emph{source cycle}  is a cycle without entrances, i.e., a cycle $c=e_1\cdots e_n$ such that, for each $e\in E^1$,  $r(e)\in c^0:=\{s(e_i):i=1,...,n\}$ implies $e\in\{e_1,...,e_n\}$.  

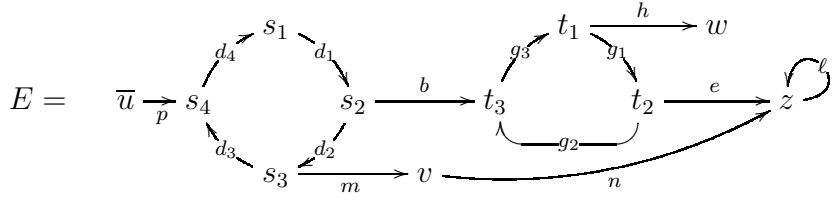
\begin{figure}
\begin{center}
\[\xymatrix@-1pc{&&&s_1\ar@/^/[dr]|{d_1}&&&&t_1\ar@/^/[dr]|{g_1}\ar[rr]^h&&w\\\
E=&\overline u\ar[r]_p&s_4\ar@/^/[ur]|{d_4}&&s_2\ar@/^/[dl]|{d_2}\ar[rr]^b&&t_3\ar@/^/[ur]|{g_3}&&t_2\ar@<-2pt> `d[l] `[ll]|{g_2} [ll]\ar[rr]^e&&z\ar@(r,u)|\ell\\
&&&s_3\ar@/^/[ul]|{d_3}\ar[rr]_m&&v\ar@/_1pc/ [rrrrru]_n}\]
\caption{Reference example of a graph {with disjoint cycles}}
\label{Figure:ref}
\end{center}
\end{figure}

Suppose $A$ is a $K$-algebra, and suppose that there is a subset  $S:=\{a_v:v\in E^0\} \sqcup \{a_e, a_{e^*}:e\in E^1\}$ of $A$ for which $\{a_v:v\in E^0\}$ is a set of orthogonal idempotents, and for which the relations analogous to  (1), (2), and (3) above are satisfied by  the elements of $S$.  Such a set $S$ is called a {\it Cuntz-Krieger} $E${\it -family} in $A$.  By the construction of $L_K(E)$, if $S$ is a Cuntz-Krieger $E$-family in $A$, then there exists a $K$-algebra homomorphism $\Phi: L_K(E) \to A$ for which $\Phi(v) = a_v$, $\Phi(e) = a_e$, and $\Phi(e^*) = a_{e^*}$ for all $v\in E^0$ and $e\in E^1$.

Because $E$ is assumed to be finite, $L_K(E)$ is a unital $K$-algebra, with multiplicative identity $1_{L_K(E)} = \sum_{v\in E^0}v$.    If $d$ is a closed path in $E$ and $p(x) = a_0 + \sum_{i=1}^n a_ix^i \in K[x]$, then $p(d)$ denotes the element 
$$ p(d) \ := \ a_0 1_{L_K(E)} + \sum_{i=1}^n a_id^i $$
of $L_K(E).$

A subset $H\subseteq E^0$ is \emph{hereditary} if, whenever $v\in H$ and there exists  $p\in {\rm Path}(E)$ for which $s(p) = v$ and $r(p)=w$, then $w \in H$. An hereditary set of vertices is \emph{saturated} if for any non-sink $v\in E^0$,  $r(s^{-1}(v))\subseteq H$ implies $v\in H$. {In the Figure~\ref{Figure:ref} the set of vertices $\{t_1,t_2,t_3,w,z\}$ is hereditary, and the set $\{v, t_1,t_2,t_3,w,z\}$ is hereditary and saturated.}
Let  $\mathcal P$ be any subset of $\text{Path}(E)$.   Let $\ell_1$, $\ell_2\in L_K(E)$ denote 
either a vertex of $E$ or the sum of distinct vertices of $E$.   We denote by $$\ell_1\mathcal P\ell_2:=\{\ell_1 p\ell_2:p\in\mathcal P\}$$  the set of  all paths in $\mathcal{P}$ which start in a vertex summand of $\ell_1$ and end in a vertex summand of $\ell_2$.  {Referring to the Figure~\ref{Figure:ref}, 
\[s_2\{g_3h,bg_3, bg_3h,bg_3g_1e\}(w+z)=\{bg_3h,bg_3g_1d\}.\]}

\smallskip

Additional general information about Leavitt path algebras may be found in \cite{AAS17}.
\medskip

 For a hereditary saturated subset $H$ of $E^0$, $E_H$ denotes the restriction of $E$ to $H$, that is, the graph with $E_H^0=H$ and $E_H^1=\{e\in E^1\mid s(e)\in H\}$. {Referring to the Figure~\ref{Figure:ref}, if $H=\{v,t_1,t_2,t_3,w,z\}$, then 
 
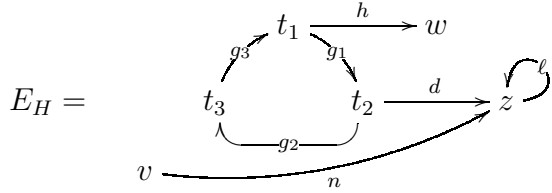
\[
 \xymatrix@-1pc{&&&t_1\ar@/^/[dr]|{g_1}\ar[rr]^h&&w\\\
E_H=&&t_3\ar@/^/[ur]|{g_3}&&t_2\ar@<-2pt> `d[l] `[ll]|{g_2} [ll]\ar[rr]^d&&z\ar@(r,u)|\ell\\
&v\ar@/_1pc/ [rrrrru]_n}
 \]}

\begin{proposition}\label{lemma:idealeprincipale}
Let $K$ be any field.   Let  $\tau$ be a source loop in  the graph $E$, with $s(\tau) = r(\tau)  = t$.   Let $H:=E^0\setminus\{t\}$, and let  $I$ be the two sided ideal of $L_K(E)$  generated by $\rho:=\sum_{u\in H}u$;  that is,  $I = L_K(E) \rho L_K(E)$.  
Then:  
\begin{enumerate}
\item  The ideal $I$, when viewed as a $K$-algebra in its own right, is isomorphic to a  Leavitt path $K$-algebra. Further, $I$ is Morita equivalent to  $L_K(E_H)$.
\item There is a lattice isomorphism between the lattice of  left $L_K(E)$-ideals contained in $I$ and the lattice of  left $L_K(E_H)$-submodules of $\rho I$.
\item Any left $L_K(E)$-ideal properly containing $I$ is a principal left ideal of $L_K(E)$.  It is  generated by
an element  of the form $p(\tau) \in L_K(E)$, where  $p(x)\in K[x]$ has $p(0)=1$. 
\end{enumerate}
%
%
%
\end{proposition}
\begin{proof}
1.)  Since $\tau$ is a source cycle, $H$ is clearly hereditary.  As $\tau$ is a loop based at $t$, $H$ is necessarily (vacuously) saturated as well.   
 Invoking \cite[Theorem 2.5.19]{AAS17} we get that  $I$ is isomorphic to a Leavitt path algebra of a (not-necessarily finite) graph ${}_HE$. (The graph ${}_HE$ is a so-called ``hedgehog graph".)
Note then that $I = L_K({}_HE)$ need not have multiplicative identity. However,  $I$ has local units when viewed as a ring (see e.g., \cite[Section 1.2]{AAS17}).   Since $\rho = \rho^2$, and $I = L_K(E) \rho L_K(E)$, clearly then $I = I \rho I$;  as well, $\rho I \rho = \rho L_K(E) \rho$.   
 Hence by  \cite[Corollary~4.3]{Ab83} we have an equivalence 
\[\xymatrix{I\text{-Mod}\ar@<1ex>[rrr]^{\rho I\otimes_I-}&&&\rho L_K(E)\rho\text{-Mod}
\ar@<1ex>[lll]^{ I\rho\otimes_{\rho L_K(E)\rho}-}}
\]
between the category $I$-Mod of unitary $I$-modules (i.e., for each $M$ in $I$-Mod, $M=IM$ holds) and the category of modules over the full corner ring $\rho L_K(E)\rho$. Since there are no paths containing $t$ 
  that both start and end in $H$, the ring $\rho L_K(E)\rho$ is isomorphic to the Leavitt path algebra $L_K(E_H)$.

\medskip

2.) Let $J$ be a left ideal of $L_K(E)$ contained in $I$. Since $I$ has local units, we have $IJ=J$ and hence $J$ belongs to $I$-Mod. Conversely, any left $I$-submodule $M$ of $I$ is a left $L_K(E)$-ideal contained in $I$, indeed
\[L_K(E)M=L_K(E)(IM)=(L_K(E)I)M=IM=M.\] 
Applying the equivalence described in point 1.), we get the indicated lattice isomorphism between the left $L_K(E)$-ideals contained in $I$ and the left $L_K(E_H)$-submodules of $ \rho I\otimes_I I \cong \rho I$.

\medskip

3.) Let $J'$ be any left ideal of $L_K(E)$ properly containing $I$.  We must prove that there exists $p(x) = 1 + \sum_{i=1}^m a_i x^i \in K[x]$ for which 
\[J'=L_K(E) p(\tau)=L_K(E)(1_{L_K(E)}+a_1\tau+\cdots+a_m\tau^m).\]
By \cite[Corollary 2.4.13]{AAS17}, $L_K(E) / I  \cong L_K(  \xymatrix{  {\bullet}^t \ar@(ul,dl)_\tau } )  \cong K[x,x^{-1}]$.
 Since $0\not=J'/I$ is then a (left) ideal of the principal ideal domain $L_K(E)/I$, it is generated as an $L_K(E)/I$-ideal by a nonzero polynomial $p(x)$ of degree $m\geq 0$ with $p(0)=1$ evaluated in  $\tau+I$. We consider the following two elements  of $L_K(E)$:
\[p_t(\tau):= t +a_1\tau+\cdots+a_m\tau^m\quad \ \ \ \ \text{and}
\]
\[p(\tau):=1_{L_K(E)}+a_1\tau+\cdots+a_m\tau^m=\rho+p_t(\tau)
.\]
We show that $J'=L_K(E)p(\tau)$. One can easily check that $J'=L_K(E)p_t(\tau)+I$, and clearly $J'$ contains $p_t(\tau)+\rho=p(\tau)$; since $p_t(\tau)= t p(\tau)$ we have
\[J'\geq L_K(E)p(\tau)\geq L_K(E)p_t(\tau).\] 
Since then $J' \leq L_K(E)p(\tau) + I$, to prove that $J'=L_K(E)p(\tau)$ we need to show that $I$ is   contained in $L_K(E)p(\tau)$. 

To do so, let $s^{-1}(t) = \{\tau, \varepsilon_1,...,\varepsilon_n\}$ be the set of edges in $E$ with source $t = s(\tau)$. Any element of $I$ is a sum of the type $\alpha \rho+\beta\rho\gamma s(\tau)$ with $\alpha,\beta,\gamma\in L_K(E)$. 
Clearly $\alpha \rho=\alpha \rho p(\tau)$ belongs to $L_K(E)p(\tau)$. Since $\tau$ is a source loop, $\rho\gamma s(\tau) =\sum_{j=1}^n k_j\rho\gamma_j\varepsilon_j^*(\tau^*)^{\ell_j}$ for suitable $k_j\in K$, $\ell_j\in\mathbb N$, and $\gamma_j\in L_K(E)$.
We prove by induction on $\ell_j$ that  $\varepsilon_j^*(\tau^*)^{\ell_j}$ belongs to $L_K(E)p(\tau)$ for any $\ell_j\geq 0$; this will then give that $I\leq L_K(E)p(\tau)$.  Note that $\varepsilon_j^* \tau = 0$ for all $1\leq j \leq n$, so  $\varepsilon_j^* p(\tau) = \varepsilon_j^*$. 

If $\ell_j=0$ then
\[ \varepsilon_j^* (\tau^*)^0 = \varepsilon_j^* = \varepsilon_j^*p(\tau)\in L_K(E)p(\tau).\]

Let $0<\ell_j\leq m=\deg p(x)$.  Then
\begin{eqnarray*}
\varepsilon_j^*(\tau^*)^{\ell_j}p(\tau) & = & \varepsilon_j^*(\tau^*)^{\ell_j}(1_{L_K(E)} +a_1\tau+\cdots+a_m\tau^m) \\ 
 & = & \varepsilon_j^*\left((\tau^*)^{\ell_j}+ a_1 (\tau^*)^{\ell_j - 1} +    \cdots+ a_{\ell_j-1} \tau^*+a_{\ell_j}t+\cdots+a_m\tau^{m-\ell_j}\right) \\
& = & \varepsilon_j^*(\tau^*)^{\ell_j}+      a_1  \varepsilon_j^*(\tau^*)^{\ell_j - 1}  +    \cdots+a_{\ell_j-1} \varepsilon_j^*(\tau^*)+\varepsilon_j^* \left(a_{\ell_j}t +\cdots+a_m\tau^{m-\ell_j}\right) \\ 
& = & \varepsilon_j^*(\tau^*)^{\ell_j}+      a_1  \varepsilon_j^*(\tau^*)^{\ell_j - 1}  +    \cdots+a_{\ell_j-1} \varepsilon_j^*(\tau^*)+(\varepsilon_j^* p(\tau))\left(a_{\ell_j}t +\cdots+a_m\tau^{m-\ell_j}\right) \\ 
 & = & \varepsilon_j^*(\tau^*)^{\ell_j} +  a_1  \varepsilon_j^*(\tau^*)^{\ell_j - 1}  +\cdots+a_{\ell_j-1}\varepsilon_j^*(\tau^*)+\varepsilon_j^*\left(a_{\ell_j}t+\cdots+a_m\tau^{m-\ell_j}\right)p(\tau).
\end{eqnarray*}
Note that in the last step we used the fact that $\tau t = t \tau = \tau$, and that polynomials in $\tau$ commute in $L_K(E)$. Therefore 
\[\varepsilon_j^*(\tau^*)^{\ell_j}=\varepsilon_j^*(\tau^*)^{\ell_j}p(\tau)- [a_1\varepsilon_j^*(\tau^*)^{\ell_j-1} + \cdots + a_{\ell_j-1}\varepsilon_j^*(\tau^*)]-\varepsilon_j^*\left(a_{\ell_j}t+\cdots+a_m\tau^{m-\ell_j}\right)p(\tau).\]
By  the inductive hypothesis, each summand in the bracketed term belongs to $L_K(E)p(\tau)$, and thus so does  $\varepsilon_j^*(\tau^*)^{\ell_j}$ . 

Finally, let $\ell_j> m$.  Then 
\begin{eqnarray*}
\varepsilon_j^*(\tau^*)^{\ell_j}p(\tau) &= & \varepsilon_j^*(\tau^*)^{\ell_j}(1+a_1\tau+\cdots+a_m\tau^m) \\
 & = & \varepsilon_j^*(\tau^*)^{\ell_j}+a_1\varepsilon_j^*(\tau^*)^{\ell_j-1}+\cdots+a_m\varepsilon_j^*(\tau^*)^{\ell_j-m}.
 \end{eqnarray*} 
Thus
\[\varepsilon_j^*(\tau^*)^{\ell_j}=\varepsilon_j^*(\tau^*)^{\ell_j}p(\tau)-[a_1\varepsilon_j^*(\tau^*)^{\ell_j-1}+\cdots + a_m\varepsilon_j^*(\tau^*)^{\ell_j-m}],\]
which, again by the inductive hypothesis, belongs to $L_K(E)p(\tau)$.
\end{proof}

The following result regarding expressions in $L_K(E)$ will be extremely useful in the sequel.

\begin{lemma}\label{lemma:tecnico2}
Let $E$ be a finite graph. 
Suppose that $\gamma_1, \gamma_2 \in {\rm Path}(E)$ have  $r(\gamma_1)= r(\gamma_2) = v$ for some $v\in E^0$.    Suppose further that neither $\gamma_1$ nor $\gamma_2$ is of the form $\gamma \delta$, where $\delta$ is a closed path having $s(\delta) = r(\delta) = v$.  
Then in $L_K(E)$, 
\[\gamma_1^*\gamma_2=\begin{cases}
    v  & \text{if }\gamma_1=\gamma_2, \\
    0  & \text{otherwise}.
\end{cases}\]
\end{lemma}

\begin{proof}
If $\gamma_1$ or $\gamma_2$ are vertices, 
then the result is easy to check.  So now assume $\gamma_1=\alpha_1\cdots\alpha_h$ and $\gamma_2=\beta_1\cdots\beta_k$, where $h,k\geq 1$, and $r(\alpha_h)=r(\beta_k) = v$. Then
\[\gamma_1^*\gamma_2=\alpha_h^*\cdots\alpha_1^*\beta_1\cdots\beta_k.\]
Assume $\gamma_1^*\gamma_2\not=0$. If $h< k$ we would have $\alpha_1=\beta_1$, ..., $\alpha_h=\beta_h$; since $r(\alpha_h)=r(\beta_h)=r(\beta_k) = v$, we get that  the path $\beta_{h+1}\cdots\beta_k$ would be a closed path with source and range vertex $v$, 
contradicting the hypothesis on $\gamma_2$.
The situation in which $h>k$ follows similarly. 
Therefore $h=k$ and $\gamma_1=\gamma_2$, and the result follows.
\end{proof}

We close this section by noting some graph-theoretic properties related to {the condition of having disjoint cycles}.  


\begin{lemma}\label{lemma:Ranga}
Let $E$ be a finite graph.   Then $E$ {has disjoint cycles} if and only if all closed paths in $E$ are powers of cycles.\footnote{This lemma has been obtained in one of our fruitful discussions with our friend and colleague Kulumani M. Rangaswamy.}
\end{lemma}
\begin{proof}
Clearly any graph in which the only closed paths are powers of cycles {has disjoint cycles}.

Conversely, let $e:=e_1\cdots e_n$ be a closed path in $E$. We proceed by induction on the length $n$ of $e$. If $n=1$, then $e$ is a loop and hence it is a cycle. Let $n>1$. Assume that $e$ is not a power of a cycle. Consider the sequence of vertices $\left(s(e_1), ..., s(e_n)\right)$. Since $e$ is not a cycle, in the sequence there are repetitions. In the set 
\[\left\{(i,j):i<j\in\{1,..., n\}, \ s(e_i)=s(e_j),  \ j-i\text{ minimal}\right\}
\]
consider the pair $(\overline i,\overline j)$ with $\overline i$ minimal. Then $e_1\cdots e_{\overline i-1}e_{\overline j}\cdots e_n$ and $e_{\overline i}\cdots e_{\overline j-1}$ are respectively a closed path $p$ and a cycle $d$ in $E$ of length $<n$. 
By the inductive hypothesis we have $p=c^m$ for a suitable cycle $c$ in $E$ and $m\geq 1$. The cycles $c$ and $d$ have in common the vertex $s(e_{\overline i})=r(e_{\overline i-1})$: therefore, {since distinct cycles are disjoint}, $c=d$. By the minimality of $\overline i$, we get $\overline i=1$ and hence
\[c=d=e_1\cdots e_{\overline j-1},\quad p=c^m=e_{\overline j}\cdots e_n.\]
Therefore $e=c^{m+1}$, a power of the cycle $c$.
\end{proof}

The set of vertices of a graph $E$ is naturally endowed with a preorder structure. Given $u,v\in E^0$, we write $v\leq u$ if there exists  $p\in {\rm Path}(E)$ with $s(p)=u$ and $r(p)=v$.

{Writing $u\equiv v$ if $u\leq v$ and $v\leq u$, one obtains an equivalence relation $\equiv$ on $E^0$. The preorder $\leq$ induces a partial order on the quotient set $E^0_\equiv$. Because  $E^0$ is assumed to be finite, $E^0_\equiv$ contains maximal elements. {Referring to the Figure~\ref{Figure:ref}, $E^0=\{\overline u,s_1,s_2,s_3,s_4,v,t_1,t_2,t_3,w,z\}$ and $E^0_\equiv=\{[\overline u],[s_1,s_2,s_3,s_4],[v],[t_1,t_2,t_3],[w],[z]\}$; $[\overline u]$ is the only maximal element in $E^0_\equiv$.
}
\begin{remark}\label{rem:rigidity}
If the graph $E$ {has disjoint cycles}, then (using Lemma \ref{lemma:Ranga}) each equivalence class in $E^0_\equiv$ is either a  single vertex,  or the set of vertices of a cycle. In particular, the maximal elements in $E^0_\equiv$ are either source vertices or the set of vertices of a source cycle.

Therefore, if $E$ in addition contains no source vertices, then the finiteness of $E$ together with {the condition of having disjoint cycles} implies that $E$ contains at least one source cycle.  
\end{remark}}
\medskip
\section{Chen simples and Pr\"ufer modules}
\label{section:ChenPrufer}

We reiterate here our standing assumption that  $E$ is a finite graph. So the Leavitt path algebra $L_K(E)$ has a multiplicative identity, namely,  $1_{L_K(E)} = \sum_{v\in E^0}v.$ 

In \cite{Ch12} Chen introduced classes of simple left modules over Leavitt path algebras. Here, we concentrate on those simples associated to sinks or to cycles: $L_K(E)w$ (for any sink $w$) and $V_{[c^\infty]}$ (for any  cycle $c$) in $E$. We will see that, although these have been viewed as distinct classes in the literature, they in fact share many properties.

A sink $w$  is an idempotent in $L_K(E)$:  we consider the left ideal $L_K(E)w$.
Since there are no ghost edges ending in $w$, any element of $L_K(E)w$ is a $K$-linear combination of paths in $\text{Path}(E)w$, i.e., of (real) paths ending in $w$. Since 
$L_K(E)=\bigoplus_{u\in E^0}L_K(E)u$ {as left $L_K(E)$-ideals}, 
we have clearly $$L_K(E)w\cong L_K(E)/L_K(E)(w-1).$$

Given a cycle $c=e_1\cdots e_n$ in $E$, the ``infinite path"  obtained by repeating $c$ infinitely many times, denoted $c^\infty$,  can be considered as an ``idempotent-like element'': one has $c^n\cdot c^\infty=c^\infty$ for each $n\geq 1$. Extending the product defined between elements of $L_K(E)$ and $c$, one can in the expected way construct the left $L_K(E)$-module $L_K(E)c^\infty$, which is usually denoted by $V_{[c^\infty]}$. Since $e^*c^\infty=0$ for each $e\not=e_1$, and $e_1^*c^\infty=e_2\cdots e_nc^\infty$, any element in $L_K(E)c^\infty$ is a $K$-linear combination of paths in $\text{Path}(E)c^\infty$, i.e., of (real)  paths ending in $s(c)$ times $c^\infty$. 
By \cite[Theorem 2.8]{AMT15} we have 
$$L_K(E)c^\infty\cong L_K(E)/L_K(E)(c-1).$$

If one thinks of a sink $w$ as a cycle of length 0, and observes that  $w^\infty=w$, the two above constructions of simple modules may be viewed as two manifestations of the same idea.   

\begin{definition}\label{cycledef}  By a {\rm cycle} in the graph $E$ we mean either a sink (i.e., a cycle of length $0$), or a {\rm proper cycle}, i.e., a cycle of length $n\geq 1$.    We unify notation by setting \[V_{[w^\infty]}:=L_K(E)w\]
for any sink $w$ in $E$.

\end{definition}

\begin{definition}\label{Pc(E)definition}
Let $E$ be a finite graph, and $c$ any cycle in $E$.  By $\mathcal P^E_c$ we denote the following subset of ${\rm Path}(E)$: 
\[  \mathcal P^E_c \ :=\begin{cases}
     \{\gamma\in {\rm Path}(E):r(\gamma)=s(c) = c\} = {\rm Path}(E)c  \\ \medskip  \hspace{.5in}  \text{if $c$ is a sink;} \\  
    \{\gamma\in {\rm Path}(E):r(\gamma)=s(c), \ {\rm and } \ \gamma\not=\gamma'c \ \ \forall\gamma'\in {\rm Path}(E)\}   \\ \hspace{.5in}  \text{if $c$ is a proper cycle}.
\end{cases}.\]
Less formally, if $c$ is a sink then $\mathcal P^E_c$ is the set of paths in $E$ ending in $s(c)=c$.  Otherwise,  if $c$ is a proper cycle,  $\mathcal P^E_c$ is the set of paths in $E$ which end at the source vertex $s(c)$ of the cycle $c$, but for which the path does not end with a complete traverse of the cycle $c$. \end{definition}

{Referring to the Figure~\ref{Figure:ref}, the set $\mathcal P^E_w$ of all paths in $E$ ending in $w$ is infinite: it contains the sink $w$ and all the left truncations of paths of the form $pd_4(d_1d_2d_3d_4)^i d_1bg_3(g_1g_2g_3)^jh$ ($ i,j\geq 0$).
As well, the sets $\mathcal P^E_{g_1g_2g_3}$,  and $\mathcal P^E_\ell$  are infinite. The former contains the vertex $t_1=s(g_1)$, $g_2g_3$, and all the left truncations of paths of the form $pd_4(d_1d_2d_3d_4)^i d_1bg_3$ ($i\geq 0$); the latter contains the vertex $z$ and all the left truncations of paths of the form
$pd_4(d_1d_2d_3d_4)^id_1d_2mn$ and $pd_4(d_1d_2d_3d_4)^id_1bg_3(g_1g_2g_3)^jg_1e$. The set $\mathcal P^E_{d_1d_2d_3d_4}$ is finite: it contains $s_1$, $d_4$, $pd_4$, $d_3d_4$, $d_2d_3d_4$.}

\begin{remark}\label{anyclosedpath}
By Lemma \ref{lemma:Ranga}, if the graph $E$ {has disjoint cycles}, then  for any cycle $c$ the paths in $\mathcal P^E_c$ do not end with the traverse of any closed path.

\end{remark}
\medskip

\begin{remark}\label{rem:base1}
It follows by \cite[\S 3.1]{Ch12} and  \cite[Corollary~1.5.15]{AAS17} that for any cycle $c$, 
 the set 
\[\mathcal P^E_c\cdot c^\infty=\{p\cdot c^\infty:p\in \mathcal P^E_c\}\] is a $K$-base for the
simple left $L_K(E)$-module $V_{[c^\infty]}$. 
\end{remark}
\medskip


{If $c$ is any cycle and $0\not=a\in K$, we denote by {$\sigma_{c,a}$} the  ``gauge'' automorphism of $L_K(E)$ associated to $c$ and $a$.  This automorphism  is defined to be the identity if $c$ is a sink.  Otherwise, if  $c= e_1 e_2 \cdots e_n$ is a proper cycle, $\sigma_{c,a}$ sends $u$ to $u$ for each $u\in E^0$, $e$ to $e$ and $e^*$ to $e^*$ for each $e\in E^1\setminus\{e_1\}$, and $e_1$ to $ae_1$ and $e_1^*$ to $ a^{-1}e_1^*$.    If $a=1$, then $\sigma_{c,1}$ is clearly the identity automorphism of $L_K(E)$.} 

{For $M$ in $L_K(E)$-Mod, we define $M^{\sigma_{c,a}}\in L_K(E)$-Mod by setting $M^{\sigma_{c,a}}=M$ as an abelian group, but with the modified left $L_K(E)$-action $$\ell\star m:=\sigma_{c,a}(\ell)m$$ for each $\ell\in L_K(E)$ and $m\in M$.
The automorphism $\sigma_{c,a}$ induces an auto-equivalence of the category $L_K(E)$-Mod given by the functor \[L_K(E)^{\sigma_{c,a}}\otimes_{L_K(E)}-: \ \ L_K(E)\text{-Mod}\to L_K(E)\text{-Mod}\]
which sends any $L_K(E)$-module $M$ to the ``twisted'' $L_K(E)$-module $M^{\sigma_{c,a}}\cong L_K(E)^{\sigma_{c,a}}\otimes M$. The twisted module $V^{\sigma_{c,a}}_{[c^\infty]}$ is a simple left $L_K(E)$-module for each $0\not=a\in K$. }

A  family of simple  $L_K(E)$-modules wider than the one presented in  \cite{Ch12}   for proper cycles 
was  obtained by Ara and Rangaswamy in \cite{AR14} as follows.  
Let  $K$ be any field, $c$ be any cycle in $E$, and $f(x) \in K[x]$ be a \emph{basic irreducible} polynomial in $K[x]$, i.e.,
a polynomial  which is irreducible in $ K[x]$, and for which $f(0)=-1$. 
Denote by $K'$ the field $K[x]/\langle f(x)\rangle$ and by $\overline x$ the element $x+\langle f(x)\rangle\in K'$. 

{\begin{definition}    We denote by $V^f_{E,[c^\infty]}$ the left $L_K(E)$-module obtained from the twisted left  $L_{K'}(E)$-module $V^{\sigma_{c,\overline x}}_{E,[c^\infty]}$ by restricting scalars  from $K'$ to $K$. 
\end{definition}  }

Here is a more precise formulation of $V^f_{E,[c^\infty]}$.  The $K$-algebra homomorphism $K\to K'$ induces a functor $\mathcal U$ from $L_{K'}(E)$-modules to $L_{K}(E)$-modules which is the right adjoint of the extension of scalars functor 
$-  \otimes_{K} K' : L_{K}(E)\text{-Mod}\to L_{K'}(E)\text{-Mod}.$  Then   $V^f_{E,[c^\infty]}$ is defined to be { $  \mathcal{U}(V^{\sigma_{c,\overline x}}_{E,[c^\infty]})$.}     

In particular, for each $0 \neq a\in K$ we may view the previously defined twisted modules in these more general terms by noting that 
{\[ V^{\sigma_{c,a}}_{E,[c^\infty]}=V^{f}_{E,[c^\infty]}\quad \text{  where }f(x)=a^{-1}x-1.\]}
We will often omit the graph $E$ when it is clear from the context, writing simply $V^f_{[c^\infty]}$. 
{If $c = w$ is a sink, then 
\[V^f_{[w^\infty]}=L_K(E)w=K\mathcal P^E_w = K{\rm Path}(E)w
\] for any basic irreducible polynomial $f(x)\in K[x]$}.

\begin{remark} Since any element of $K'$ is a $K$-linear combination of elements of the form $\overline x^h \  (0\leq h < deg(f(x)))$, by Remark~\ref{rem:base1} for any cycle $c$,
 the set 
\[\mathcal P^E_c\cdot\{\overline x^i:0\leq i<\deg f(x)\}\cdot c^\infty=\{p\cdot \overline x^i\cdot c^\infty:p\in \mathcal P^E_c, \  0\leq i<\deg f(x)\}\] is a $K$-base for the
simple left $L_K(E)$-module $V^f_{[c^\infty]}$. 
\end{remark}
\medskip

\begin{definition} Let $c$ be any sink and $f(x)$ any basic irreducible polynomial in $K[x]$.    Following Ara and Rangaswamy \cite{AR14},  we call any simple left $L_K(E)$-module of the form $V^f_{[c^\infty]}$
a {\rm Chen simple module}.  
\end{definition}

\begin{remark}  If $c_1$ is a cycle of length $\geq 2$ obtained by ``rotating'' the cycle $c_2$, then the simple modules $V^f_{E,[c_1^\infty]}$ and $V^f_{E,[c_2^\infty]}$ are isomorphic. If $c$ is a fixed sink, then  $V^f_{E,[c^\infty]}=V_{E,[c^\infty]}$, for any basic irreducible polynomial $f(x)\in K[x]$. Aside from  these cases, the Chen simple modules are pairwise non-isomorphic.  (See \cite[Theorem 6.2]{Ch12} and \cite[Proposition 3.8]{AR14}.)
\end{remark}

\medskip
 {Referring to the Figure~\ref{Figure:ref}, a complete list of 
the simple left $L_K(E)$-modules up to isomorphisms is given by:    $V_{[w^\infty]}=L_K(E)w$,  $V^f_{[(d_1d_2d_3d_4)^\infty]}$, $V^f_{[(g_1g_2g_3)^\infty]}$, and $V^f_{[\ell^\infty]}$, for any basic irreducible polynomial $f(x)\in K[x]$.}


 \medskip

In addition to the Chen simple modules, there is another important (and related) construction of $L_K(E)$-modules which will be relevant in the sequel.  Let $E$ be a finite graph, and $c$ any 
cycle in $E$.    
  In \cite{AMT19} the three authors  constructed  the \emph{Pr\"ufer  $L_K(E)$-module} $U_{E,c-1}$ associated to any proper cycle $c$ in $E$.  
  The module  $U_{E,c-1}$   is  an infinite length uniserial artinian $L_K(E)$-module with all composition factors isomorphic to $V_{[c^\infty]}$. The module $U_{E,c-1}$ is the direct limit 
  $$U_{E,c-1} := \  \varinjlim_{i,j\geq 1} \ \{L_K(E)/L_K(E)(c-1)^i, \  \psi_{E,i,j}\}$$
  of the factor modules  $L_K(E)/L_K(E)(c-1)^i$,  with respect to the morphisms
\[\psi_{E,i,j}:L_K(E)/L_K(E)(c-1)^i\to L_K(E)/L_K(E)(c-1)^j, \ \mbox{defined by setting}\]
{\[1+L_K(E)(c-1)^i \ \ \mapsto \ \  \frac{(c-1)^j}{(c-1)^i}+L_K(E)(c-1)^j\quad\forall 1\leq i<j.\]}
Denoting by 
$$\psi_{E,i}: \  L_K(E)/L_K(E)(c-1)^i\to U_{E,c-1}$$
 for each $i\geq 1$  the induced monomorphism, the Pr\"ufer module $U_{E,c-1}$ is generated as an $L_K(E)$-module by the elements
\[\alpha_{c,i}:= \ \psi_{E,i}(1+L_K(E)(c-1)^i),\quad i\geq 1.\]
{If $c=w$ is a sink, then $\alpha_{w,i}=(-1)^{i+1}\alpha_{w,1}$.  The  Pr\"ufer module corresponding to a sink then becomes 
 \[U_{E,w-1}=V_{[w^\infty]}=L_K(E)w=K\mathcal P^E_w.\]}
%


Clearly, for any cycle $c$,  $L_K(E)\alpha_{c,1}$ is isomorphic to the simple module $V_{[c^\infty]}$, and hence it is the socle of the module $U_{E,c-1}$ (it is equal to $V_{[c^\infty]}$ if $c$ is a sink). In particular $V_{[c^\infty]}$ is essential in $U_{E,c-1}$.
In general we omit reference to the graph $E$ and write simply $U_{c-1}$.

We now provide some additional information about the Pr\"{u}fer $L_K(E)$-module  $U_{c-1}$ for any cycle $c$ of $E$.  
{Let $m \geq 1$. Setting $\alpha_{c,0}=0$, we have  
\[(c-1)\alpha_{c,m}=\begin{cases}
  0    & \text{if $c$ is a sink}, \\
 \alpha_{c,m-1}     & \text{if $c$ is a proper cycle}.
\end{cases}\]}

   If $u\in E^0\setminus\{s(c)\}$,   then easily one shows that   $u=(-1)^mu(c-1)^m$ for each $m\geq 1$, so we get   $u\alpha_{c,m}=0$.   (As a result,   $s(c)\alpha_{c,m}=\alpha_{c,m}$.)

  If $c$ has length 0 and $e\in E^1$, or $c$ is a proper cycle and $e\in E^1\setminus\{e_1\}$, it is also easily shown that  $e^*=(-1)^me^*(c-1)^m$ for each $m\geq 1$.  So  $e^*\alpha_{c,m}=0$ for all $m\geq 1$ in this situation as well.   
    
 Finally, if  $c$ is a proper cycle $c = e_1 e_2 \cdots e_n$, then for each $m\geq 1$, multiplying each side of the equation $(c-1)\alpha_{c,m} = \alpha_{c,m-1}$ by $e_1^*$ and  rearranging terms, we get that   $e_1^*\alpha_{c,m}= e_2\cdots e_n \alpha_{c,m} -e_1^*\alpha_{c,m-1}.$  (We interpret $e_2 \cdots e_n$ as $t=s(c)$ if $n=1.$)   Continuing in this way, we iteratively get 
\[e_1^*\alpha_{c,m}=e_2\cdots e_n\left(\sum_{\ell=0}^{m-1}(-1)^\ell \alpha_{c,m-\ell}\right).\]

Therefore, by the three previous paragraphs, any element of $U_{c-1}$ can be written in the form
\[\sum_{j=1}^m\sum_{j\in F} k_{j,\gamma} \gamma\alpha_{c,j}\]
for  suitable $k_{j,\gamma}\in K$ and finite subset 
$F\subseteq \mathcal P^E_c$.
Less formally, we have established that each element of $U_{c-1}$ can be written as an $L_K(E)$-linear combination of the $\alpha_{c,m}$ in such a way that the $L_K(E)$-coefficients  do not 
involve any ghost edges.   This observation will play a key role later on.  

Analogous to the construction of the simple module $V^f_{[c^\infty]}$ associated to any proper cycle $c$ and any basic irreducible $f(x)\in K[x]$, we have
\begin{definition}\label{Uf(c)def}   Let $E$ be any finite graph, and $K$ any field.  Let $c$ be any cycle, $f(x)$ be any basic  irreducible polynomial in $K[x]$, and $\overline x:=x+\langle f(x)\rangle\in K'=K[x]/\langle f(x)\rangle$.  We define 
\[U_{f(c)}:=\mathcal U(U^{\sigma_{c,\overline x}}_{c-1}),\]
the left $L_K(E)$-module obtained from the twisted $L_{K'}(E)$-module {$U^{\sigma_{c,\overline x}}_{E, c-1}$} via the functor $\mathcal{U}$ which restricts scalars to $K$.
%
\end{definition}

{Observe that if $c=w$ is a sink, then Definition \ref{Uf(c)def} reduces to   
\[U_{f(w)}=U_{w-1}=V_{[w^\infty]}=L_K(E)w=K\mathcal P^E_w.\]}

\begin{remark}\label{rem:multbyoverx}
Since for each $0\not=a\in K$, any automorphism of $L_K(E)$ induces an autoequivalence of $L_K(E)$-Mod, it is clear that for any proper cycle $c=e_1\cdots e_n$ the twisted $L_{K}(E)$-module $U^{\sigma_{c,a}}_{E, c-1}$ has the same submodule structure as $U_{E, c-1}$. Precisely, $U^{\sigma_{c,a}}_{E, c-1}$ is an infinite length uniserial artinian left $L_K(E)$-module with composition factors isomorphic to $V^{\sigma_{c,a}}_{[c^\infty]}$. {Moreover, setting
\[\alpha^a_{c,m}:=1_{L_K(E)}\otimes_{L_K(E)}\alpha_{c,m}\quad\text{in}\quad L_K(E)^{\sigma_{c,a}}\otimes_{L_K(E)}U_{E, c-1}
\cong U^{\sigma_{c,a}}_{E, c-1}\]}
we have 
\begin{eqnarray*}
(a^{-1}c-1)\alpha^a_{c,m} & = & \left((a^{-1}c-1)\star 1_{L_K(E)}\right)\otimes_{L_K(E)}\alpha_{c,m} \\
 & = & \sigma_{c,a}(a^{-1}c-1)\otimes_{L_K(E)}\alpha_{c,m} \\
 & = &(c-1) \otimes_{L_K(E)}\alpha_{c,m}\\
  & = &1_{L_K(E)}\otimes_{L_K(E)} (c-1)\alpha_{c,m} \\
  & = & \alpha^a_{c,m-1}.
\end{eqnarray*}
If $f(x)\in K[x]$ is basic irreducible and $K'=K[x]/\langle f(x)\rangle$, taking $a$ equal to the invertible element $\overline x=x+\langle f(x)\rangle$ of $K'$  we get that 
\[(\overline x^{-1}c-1)\alpha^{\overline x}_{c,m}=\alpha^{\overline x}_{c,m-1}\]
in the twisted left $L_{K'}(E)$-module $U^{\sigma_{c,\overline x}}_{E, c-1}$. Multiplying by $\overline x$ on both sides one gets $\overline x  \alpha_{c,m}^{\overline x}=c \alpha_{c,m}^{\overline x}-\overline x  \alpha_{c,m-1}^{\overline x}$, which iteratively gives
\[\overline x\alpha^{\overline x}_{c,m}=c\alpha^{\overline x}_{c,m}-c\alpha^{\overline x}_{c,m-1}+\cdots+(-1)^{m-1}c\alpha^{\overline x}_{c,1}.\]
In such a way, multiplication of $\alpha^{\overline x}_{c,m}$  by $\overline x$ yields a linear combination of $\alpha^{\overline x}_{c,j}$'s having coefficients in $L_K(E)$.
{Moreover for each edge $e\not=e_1\in E^1$ we have
\[e^*\overline x^h\alpha^{\overline x}_{c,m}=\overline x^h e^*\alpha^{\overline x}_{c,m}=\overline x^h (1_{L_K(E)}\otimes_{L_K(E)}e^*\alpha_{c,m})=0, 
\]
while when $e=e_1$ we have 
\begin{eqnarray*}
e_1^*\overline x^h\alpha^{\overline x}_{c,m} & =& \overline x^h e_1^*\alpha^{\overline x}_{c,m} = \overline x^h e_1^*(1_{L_K(E)}\otimes_{L_K(E)}\alpha_{c,m})  \\
& = & \overline x^h (1_{L_K(E)}\otimes_{L_K(E)}e_1^*\alpha_{c,m}) \\ 
& =& \overline x^h (1_{L_K(E)}\otimes_{L_K(E)}e_2\cdots e_n\sum_{\ell=0}^{m-1}(-1)^\ell \alpha_{c,m-\ell}) \\
& = & \sum_{\ell=0}^{m-1}(-1)^\ell e_2\cdots e_n\overline x^h\alpha^{\overline x}_{c,m-\ell}.
\end{eqnarray*}
}


\end{remark}
\medskip

The discussion over the previous paragraphs has established the following.

\begin{proposition}\label{descriptionofUasdirectsum}
For any sink $w$,  any proper cycle $c$ and any basic irreducible $f(x)\in K[x]$ we have the following equalities of $K$-vector spaces
\[U_{w-1}=K\mathcal P^E_w\alpha_{w,1}=K\mathcal P^E_w\quad\text{and}\quad U_{f(c)}=\bigoplus_{j\geq 1}\sum_{h=0}^{\deg f(x)-1}K\mathcal P^E_c\overline x^h\alpha^{\overline x}_{c,j}.\]
{ Specifically, any element of $U_{f(c)}$ is a finite sum of elements of the form 
$\{k\gamma\overline x^h\alpha^{\overline x}_{c,j} \ | \ k\in K, \  \gamma\in \mathcal P^E_c, \  h\geq 0, \  j \geq 1 \}$. 
}
\end{proposition}
\medskip


\medskip

As mentioned previously, the cases of simple $L_K(E)$-modules corresponding to sinks and those corresponding to proper cycles have historically been treated separately in the literature.  We believe that there is enough commonality to these two types of simples that merits  treating them as two cases of the same construction.   
Indeed, this is the case in all of our subsequent results.  Thus, as we already said in Definition~\ref{cycledef}, when we use the word \emph{cycle} we will mean either a sink, or a proper cycle.
Admittedly, however,  there is a small price of additional notation to pay   in order to achieve this single approach.     
\begin{definition}
We set
\[\mathbb I_c:=\begin{cases}
    \mathbb N_{\geq 1}  & \text{if $c$ is a proper cycle } \\
    \{1\}  & \text{if $c=w$ is a sink}
\end{cases}\]
\[ {\cdeg f(x):=\begin{cases}
     \deg f(x) & \text{if $c$ is a proper cycle }, \\
     1 & \text{if $c=w$ is a sink}.
\end{cases}}\]
Moreover, if $c=w$ is a sink, then $\alpha^{\overline x}_{c,1}:=w$.
\end{definition}


\begin{lemma}\label{lemma:base}
Let $c$ be any cycle in the graph $E$, and let $f(x) \in K[x]$ be basic irreducible.  Assume the graph $E$ {has disjoint cycles}. Then in $U_{f(c)}$ the set 
$$\{\gamma\overline x^h\alpha^{\overline x}_{c,j} \ | \ \gamma\in \mathcal P^E_c, \ 0\leq h<\cdeg f(x), \  j\in\mathbb I_c\}$$
is  $K$-linearly independent.
\end{lemma}
\begin{proof}
Assume
\[\sum_{j\in F_1}\sum_{h=0}^{\cdeg f(x)-1}\sum_{\gamma\in F_2}k_{j,\gamma,h}\gamma\overline x^h\alpha^{\overline x}_{c,j}=0\]
for suitable finite subsets  $F_1$ of $\mathbb I_c$ and $F_2$ of $\mathcal P^E_c$. {Since the cycles are disjoint,} we can apply Lemma~\ref{lemma:Ranga} to get that the elements in $\mathcal P^E_c$ do not end with any closed path. Therefore, multiplying on the left by any $\overline\gamma^*$ (with $\overline \gamma \in \mathcal P^E_c$), by Lemma~\ref{lemma:tecnico2} we get
\[\sum_{j\in F_1}\left(\sum_{h=0}^{\cdeg f(x)-1}k_{j,\overline\gamma,h}\overline x^h\right)\alpha^{\overline x}_{c,j} \ = 0.\]
Since the set  $\{\alpha^{\overline x}_{c,j} \ | \ j \in F_1\}$ is linearly independent over the field extension $K'=K[x]/\langle f(x)\rangle$ of $K$,   
we have that 
$\sum_{h=0}^{\cdeg f(x)-1}k_{j,\overline\gamma,h}\overline x^h = 0$  in $K'$ for each $1\leq j\leq m$. Since $\{\overline x^h \ | \ 0\leq h<\cdeg f(x)\}$ is $K$-linearly independent, we get $k_{j,\overline\gamma,h}=0$ for any $j\in F_1$ and $0\leq h<\cdeg f(x)$. Because $\overline \gamma$ was arbitrary, we are done.   
\end{proof}

As expected, the wider class of modules $U_{f(c)}$  satisfy properties similar to the specific modules $U_{c-1}$.   
\begin{proposition}\label{Ufcfinflengthuniserial}
Let $E$ be a finite graph, $K$ any field, $c$ a proper cycle in $E$, and $f(x)=xg(x)-1$ a basic irreducible polynomial in $K[x]$.   The  $L_K(E)$-module $U_{f(c)}$ is an infinite length uniserial artinian module with all composition factors isomorphic to $V^f_{[c^\infty]}$.
\end{proposition}
\begin{proof}
As sets, we have $U_{f(c)}=U_{g(\overline x)c-1}$.
It is sufficient to prove  that the lattice of $L_K(E)$-submodules of $U_{f(c)}$ is equal to the lattice of $L_{K'}(E)$-submodules of $U_{g(\overline x)c-1}$. Clearly any $L_{K'}(E)$-submodule of $U_{g(\overline x)c-1}$ is also a $L_K(E)$-submodule of $U_{f(c)}$. Consider now an $L_K(E)$-submodule $M$ of $U_{f(c)}$. Since, as observed in Remark~\ref{rem:multbyoverx},
\[\overline x\alpha^{\overline x}_{c,j}=c\alpha^{\overline x}_{c,j}-c\alpha^{\overline x}_{c,j-1}+\cdots+(-1)^{j-1}c\alpha^{\overline x}_{c,1},\]
we have that $\overline x m$ belongs to $M$ for each $m\in M$. Then $M$ is also a $L_{K'}(E)$-submodule of $U_{g(\overline x)c-1}$.
\end{proof}

{ We end the section with an observation about the left-module  action of  $L_K(E)$ on $U_{f(c)}$. }
 
\begin{remark}\label{equalonBremark}
{\rm Let $c$ be either a sink in $E$, or $c = e_1 e_2 \cdots e_n$  a cycle in $E$.  Let $e\in E^1$ and $\gamma \in \mathcal P_c^E$ for which $r(e) = s(\gamma)$.   In general there are two possible interpretations of an expression of the form  $e\gamma$.  One is simply the concatenation of the paths $e$ and $\gamma$.  The other is as $e\cdot \gamma$, where we view $e\in L_K(E)$,  $\gamma \in U_{f(c)}$, and $\cdot$ denotes the left $L_K(E)$-action on $U_{f(c)}$.    

 If for example $e = e_1$ and $\gamma = e_2 \cdots e_n$ then $\gamma \in \mathcal P_c^E$, but the concatenation $e\gamma$ is not in $\mathcal P_c^E$.    So in this situation the concatenation $e\gamma$  is not an element of the form indicated in Proposition \ref{descriptionofUasdirectsum}, and the two interpretations of $e\gamma$ differ.    

However,  consideration of the lengths of  paths shows that the situation described in the previous paragraph is the only configuration one needs to consider in order to avoid the situation in which the concatenation of an edge $e$ with a path $\gamma \in \mathcal P_c^E$ having $r(e) = s(\gamma)$ yields a path which is not in $\mathcal P_c^E$.   Consequently, suppose  for example that $\gamma \in \mathcal P_c^E$ for which $|\gamma| \geq |c|$, and let $e\in E^1$ with $r(e) = s(\gamma)$.    Then the concatenation path $e \gamma$ is necessarily in $\mathcal P_c^E$.  
Thus for an element $y$ of $U_{f(c)}$ of the form  $y = k \gamma  \overline x^h  \alpha_{c,j}$ with $k\in K$, $|\gamma| \geq |c|$ and $r(e) = s(\gamma)$,  the  element $e\cdot y$ of $U_{f(c)}$ is equal to 
the element $k(e\gamma)\alpha_{c,j}$ of $U_{f(c)}$.  

A similar observation holds for $e^*$ where $e\in E^1$.   Specifically, suppose $\gamma \in \mathcal P_c^E$ has $\gamma = e \gamma^\prime$ for some $\gamma^\prime \in {\rm Path}(E)$.  
Then necessarily $\gamma^\prime \in \mathcal P_c^E$.   So the two possible interpretations of the expressions  $e^* \gamma$  (as either the element $\gamma^\prime$ of $\mathcal P_c^E$, or as the element $e^* \cdot \gamma$ of $U_{f(c)}$) coincide.   We note that the only elements of  the form $k \gamma  \overline x^h  \alpha_{c,j}$ of $U_{f(c)}$ for which $e^* k \gamma \overline x^h \alpha_{c,j} \neq 0$ and $\gamma \neq e \gamma^\prime$ must have $\gamma = v = s(c)$.   In this case, the iterative version of the expression for $e_1^*  \overline x^h \alpha_{c,m}$ 
given in Remark \ref{rem:multbyoverx}
yields a sum of elements of the indicated form, where the lengths of the paths appearing are each equal to  $ |c| - 1$.   }
\end{remark}
\medskip
 %
 %

\section{Formal power series built from $E$}\label{fps} 

In this section we introduce the key construction by which we will produce the injective envelopes of Chen simple modules.

For any cycle $c$, the set $\mathcal P^E_c\subseteq {\rm Path}(E)$ is a $K$-linearly independent set in $L_K(E)$ by \cite[Corollary~1.5.15]{AAS17}.

\begin{definition}
We denote by $K[[\mathcal P^E_c]]$ the $K$-vector space of all mappings from $\mathcal P^E_c$ to $K$. Any element $\mathfrak p$ in $K[[\mathcal P^E_c]]$ can be represented as a  ``\emph{$\mathcal P^E_c$-formal series}''
\[ \mathfrak p = \sum_{\mu\in \mathcal P^E_c}\mathfrak p(\mu)\cdot\mu.\]
\end{definition}
If $ \mathfrak p(\mu)\not=0$ only for a finite number of $\mu\in\mathcal P^E_c$, then $\mathfrak p$ belongs to the $K$ vector space $K\mathcal P^E_c$ generated by $\mathcal P^E_c$. Otherwise $\mathfrak p$ is called a \emph{proper} $\mathcal P^E_c$-formal series.

\begin{definition}
Let  $E$ be any finite graph, $c$ a  cycle  in $E$, $K$ any field,  and $f(x)$ a basic irreducible polynomial in $K[x]$.  
We  define the $K$-vector space $\widehat U_{f(c)}$  by setting 
 \[\widehat U_{f(c)} \ := \ \bigoplus_{j\in\mathbb I_c}\sum_{h=0}^{\cdeg f(x)-1}K[[\mathcal P^E_c]]\overline x^h\alpha^{\overline x}_{c,j}.
 \]
\end{definition}

\medskip

\noindent
Clearly $\widehat U_{f(c)}$ 
is an enlarged version of the Pr\"ufer module \[U_{f(c)}=\bigoplus_{j\in\mathbb I_c}\sum_{h=0}^{\cdeg f(x)-1}K\mathcal P^E_c\overline x^h\alpha^{\overline x}_{c,j}
.\]


\begin{proposition}
The $K$-vector space 
$\widehat U_{f(c)}$
is a left $L_K(E)$-module.
\end{proposition}
\begin{proof}
The set of paths $\mathcal P^E_c$ is a disjoint union
\[\{\mu\in \mathcal P^E_c:|\mu|\leq |c|\}  \ \ \sqcup  \  \ \{\mu\in \mathcal P^E_c:|\mu|= |c|+1\} \ \  \sqcup \ \ 
\{\mu\in \mathcal P^E_c:|\mu| \geq |c|+2\},\]
which we denote respectively by  $\mathcal P^E_{\leq|c|}$, $\mathcal P^E_{=|c|+1}$, and $\mathcal P^E_{ \geq |c|+2}$.   This disjoint union induces a $K$-vector space decomposition   $\widehat U_{f(c)} = A  \oplus B \oplus C$, where  
\[A:=\bigoplus_{j\in\mathbb I_c}\sum_{h=0}^{\cdeg f(x)-1}K\mathcal P^E_{\leq|c|}\overline x^h\alpha^{\overline x}_{c,j},\ B:=\bigoplus_{j\in\mathbb I_c}\sum_{h=0}^{\cdeg f(x)-1}K\mathcal P^E_{=|c|+1}\overline x^h\alpha^{\overline x}_{c,j}, \]
\[\ \ \ \text{ and} \ \ \ \  C:=\bigoplus_{j\in\mathbb I_c}\sum_{h=0}^{\cdeg f(x)-1}K\mathcal P^E_{\geq |c|+2}\overline x^h\alpha^{\overline x}_{c,j}.
\]
We note that $A \oplus B \subseteq U_{f(c)}$.   Additionally, for $\mu \in \mathcal P^E_{=|c|+1} \sqcup  \mathcal P^E_{ \geq |c|+2} \ := \ \mathcal P^E_{ \geq |c|+1}$ 
{and $e\in E_1$ we have either $e\mu=0$ or, having $r(e)=s(\mu)$, $e\mu\in\mathcal P^E_c$}, more specifically, $e\mu \in \mathcal P^E_{ \geq |c|+2}$.  

So each element of $\widehat U_{f(c)}$ is a finite sum of  elements of the form
\[\sum_{\mu_{h,j}\in\mathcal P_{c}^E}\mathfrak p(\mu_{h,j}) \mu_{h,j}\overline x^h\alpha^{\overline x}_{c,j} \in A \oplus B \oplus C, \ \  \mbox{written as} \] \[\sum_{\mu_{h,j}\in\mathcal P_{\leq |c|}^E}\mathfrak p(\mu_{h,j}) \mu_{h,j}\overline x^h\alpha^{\overline x}_{c,j}+
\sum_{\mu_{h,j}\in\mathcal P_{=|c|+1}^E}\mathfrak p(\mu_{h,j}) \mu_{h,j}\overline x^h\alpha^{\overline x}_{c,j}+\sum_{\mu_{h,j}\in\mathcal P_{\geq |c|+2}^E}\mathfrak p(\mu_{h,j}) \mu_{h,j}\overline x^h\alpha^{\overline x}_{c,j}\]
where $0\leq h<\deg f(x)$, $j\in F_{h,c}$ for some finite subset  $F_{h,c}$  of $\mathbb I_c$,  and $\mathfrak p(\mu_{h,j}) \in K$ for all $\mu_{h,j}\in\mathcal P_c^E$.

For $v \in E^0$ we define 
$P_v: \widehat U_{f(c)} \to \widehat U_{f(c)} $ by setting 
$$P_v  \left(  \sum_{\mu_{h,j}\in\mathcal P_{c}^E}\mathfrak p(\mu_{h,j}) \mu_{h,j}\overline x^h\alpha^{\overline x}_{c,j}\right)   \ \ = \ \ \sum_{\substack{\mu_{h,j}\in\mathcal P_{c}^E \\ s(\mu_{h,j}) = v}}\mathfrak p(\mu_{h,j}) \mu_{h,j}\overline x^h\alpha^{\overline x}_{c,j} .$$
%
%

Let $\cdot$ denote the left $L_K(E)$-action on $U_{f(c)}$, and  let $\iota$ denote the inclusion map from $U_{f(c)}$ to $\widehat U_{f(c)}$.  

\smallskip

For $e\in E^1$  we define $S_e: \widehat U_{f(c)} \to \widehat U_{f(c)}$ as follows.   We define two maps:   first, 
$S_{e, A\oplus B}: A\oplus B \to \widehat U_{f(c)}$ by setting 
$$ S_{e, A\oplus B}(a + b) = \iota ( e \cdot (a+b)),  $$ 
and second, $ S_{e,B\oplus C} :  B\oplus C \to \widehat U_{f(c)}$ by setting  
\[S_{e,B\oplus C} \left(\sum_{\mu_{h,j}\in\mathcal P_{\geq |c|+1}^E}\mathfrak p(\mu_{h,j}) \mu_{h,j}\overline x^h\alpha^{\overline x}_{c,j}\right)=\sum_{\substack{\mu_{h,j}\in\mathcal P^E_{\geq |c|+1} \\ r(e) = s(\mu_{h,j})}}\mathfrak p(\mu_{h,j}) e \mu_{h,j}\overline x^h\alpha^{\overline x}_{c,j}.\]
By Remark \ref{equalonBremark},   $S_{e,A\oplus B}$ and $S_{e,B\oplus C}$ coincide on $B$.  So taken together these two maps yield $S_e: \widehat U_{f(c)} \to \widehat U_{f(c)}$.  We note that  $S_e(A) \subseteq A \oplus B$ and $S_e(B\oplus C) \subseteq C$.  

\smallskip

For $f\in E^1$ we define  $S_{f^*}: \widehat U_{f(c)} \to \widehat U_{f(c)}$ as follows.   Let $\mathcal P_c^E(f) $ denote $\{\mu \in \mathcal P_c^E \ | \ \mu = f \mu^\prime\}$ for some $\mu^\prime \in {\rm Path}(E)$; note that   in this case necessarily $\mu^\prime \in \mathcal P_c^E$.    We define two maps:   first,   $S_{f^*, A\oplus B}: A\oplus B \to \widehat U_{f(c)}$ by setting 
$$ S_{f^*, A\oplus B}( a + b) = \iota ( f^* \cdot (a+b)), $$
and second,    $S_{f^*, B\oplus C}: B\oplus C \to \widehat U_{f(c)}$  by setting 
     \[ S_{f^*,B\oplus C}\left(\sum_{\mu_{h,j}\in\mathcal P_{\geq |c|+1}^E}\mathfrak p(\mu_{h,j}) \mu_{h,j}\overline x^h\alpha^{\overline x}_{c,j}\right)=
\sum_{\substack{\mu_{h,j}\in\mathcal P_{\geq |c|+1}^E(f) \\ \mu_{h,j} = f \mu^\prime_{h,j}}}\mathfrak p(\mu_{h,j})  \mu^\prime_{h,j}\overline x^h\alpha^{\overline x}_{c,j}.\]
Again using Remark \ref{equalonBremark},    $S_{f^*, A\oplus B}$ and $S_{f^*,B\oplus C}$ coincide on $B$.  So taken together these two maps yield $S_{f^*}: \widehat U_{f(c)} \to \widehat U_{f(c)}$.  We note that $S_{f^*}(A\oplus B) \subseteq A$.


\smallskip

{\bf Claim}:   The subset $\{P_v, S_e, S_{e^*} \ | \ v\in E^0, e \in E^1\}$  forms a Cuntz Krieger $E$-family in ${\rm End}_K(\widehat U_{f(c)})$.

  {\bf Proof of Claim}:   We suppress the composition operator $\circ$ in the following equations;  $\delta$ denotes the Kronecker delta.    By definition it must be checked that:
\begin{enumerate}
\item  $P_v P_w=\delta_{v,w} P_v$ for all $v,w\in E^0$; 
\item  $P_{s(e)} S_e = S_e = S_e P_{r(e)}$ for all $e\in E^1$;
\item  $P_{r(e)} S_{e^*} = S_{e^*} =  S_{e^*} P_{s(e)}$ for all $e\in E^1$;
\item  $S_{f^*}S_e = \delta_{e,f}P_{r(e)}$ for all $e,f \in E^1$;  \ and 
\item   $\sum_{e \in s^{-1}(v)} S_e S_{e^*} =
 P_v$ for all non-sink $v \in E^0$.
\end{enumerate}

Statements (1), (2), and (3) are straightforward to check.  We give now the key details of the verification of statement (4).   We use the decomposition of $\widehat{U}_{f(c)}$ given at the start of the proof of the proposition.     

\medskip

Proof of (4). \    Let $a\in A$.  Recall that $S_e(a) \in A\oplus B$.    Then 
$$S_{f^*}S_e(a) = S_{f^*} (e \cdot a) = S_{f^*, A\oplus B}(e \cdot a) = f^* \cdot (e\cdot a) $$
$$= (f^* e) \cdot a = \delta_{e,f} r(e) \cdot a = \delta_{e,f} P_{r(e)}(a),$$ so that $S_{f^*}S_e = \delta_{e,f}P_{r(e)}$ on $A$. 

    Now let $y\in B\oplus C$.  Then $y$ is a finite sum of terms of the form 
$ \sum_{\mu_{h,j}\in\mathcal P_{ \geq |c|+1}^E}\mathfrak p(\mu_{h,j}) \mu_{h,j}\overline x^h\alpha^{\overline x}_{c,j} $, 
where $0\leq h<\deg f(x)$, $j\in F_{h,c}$ where $F_{h,c}$ is some finite subset of $\mathbb I_c$,  and $\mathfrak p(\mu_{h,j}) \in K$ for all $\mu_{h,j}\in\mathcal P_{\geq |c|+1}^E$.   Recall that $S_e(y) \in C$.   Then
\[S_e(y) = S_{e, B\oplus C}(y) = S_{e,B\oplus C} \left(\sum_{\mu_{h,j}\in\mathcal P^E_{\geq |c|+1} }\mathfrak p(\mu_{h,j}) \mu_{h,j}\overline x^h\alpha^{\overline x}_{c,j}\right) \]
\[=\sum_{\substack{\mu_{h,j}\in\mathcal P^E_{\geq |c|+1} \\ r(e) = s(\mu_{h,j})}}\mathfrak p(\mu_{h,j}) e \mu_{h,j}\overline x^h\alpha^{\overline x}_{c,j}.\]

If $f\neq e$, then  $\{ e \mu_{h,j} \ | \ r(e) = s(\mu_{h,j})\} \ \cap \  \mathcal P_{\geq |c|+1}^E(f) \ = \emptyset$ by definition, so that  
$$S_{f^*}(S_e(y)) =  S_{f^*, B\oplus C} \left( \sum_{\substack{\mu_{h,j}\in\mathcal P^E_{\geq |c|+1} \\ r(e) = s(\mu_{h,j})}}\mathfrak p(\mu_{h,j}) e \mu_{h,j}\overline x^h\alpha^{\overline x}_{c,j}\right)  = 0.$$  
On the other hand, since $S_e(y) \in C$, we have 
$$S_{e^*}(S_e(y)) =  S_{e^*, B\oplus C} \left( \sum_{\substack{\mu_{h,j}\in\mathcal P^E_{\geq |c|+1} \\ r(e) = s(\mu_{h,j})}}\mathfrak p(\mu_{h,j}) e \mu_{h,j}\overline x^h\alpha^{\overline x}_{c,j}\right)  $$
$$=  \sum_{\substack{\mu_{h,j}\in\mathcal P^E_{\geq |c|+1} \\ r(e) = s(\mu_{h,j})}}\mathfrak p(\mu_{h,j})  \mu_{h,j}\overline x^h\alpha^{\overline x}_{c,j} = P_{r(e)}(y).$$
Thus we have that $S_{f^*}S_e = \delta_{e,f}P_{r(e)}$ on $B \oplus C$ as well.    
We conclude that $S_{f^*}S_e = \delta_{e,f}P_{r(e)}$ on all of $\widehat{U}_{f(c)}$, as desired.

\medskip

Proof of (5).  \ This follows in a manner similar to the proof of (4), and is left to the reader.     The key observation here is that   for each non-sink $v\in E^0$ we have $\{\mu \in \mathcal P_c^E \  | \  s(\mu) = v\} = \bigsqcup_{e\in s^{-1}(v)} \mathcal P_c^E(e)$.  (Less formally:  every path in $E$ having source vertex $v$ must have as its first edge one of the finitely many elements of $s^{-1}(v)$.)

\smallskip

With the Claim  established, we invoke the universal property of $L_K(E)$ to conclude there exists a $K$-algebra homomorphism
$$\Phi : L_K(E) \to {\rm End}_K(\widehat U_{f(c)})$$
for which $\Phi(v) = P_v, \Phi(e) = S_e, $ and $\Phi(e^*) = S_{e^*}$ for all $v\in E^0, e \in E^1$.  

Since every path in $E$ has source vertex in $E^0$, it is clear that  $1_{{\rm End}_K(\widehat U_{f(c)})} = \sum_{v\in E^0}P_v$.    Thus $\Phi(1_{L_K(E)}) =  \Phi(\sum_{v\in E^0} v) = 1_{{\rm End}_K(\widehat{U}_{f(c)})}$.  

\medskip

With the above discussion and computations in hand, by a standard ring-theory argument we conclude that $\widehat U_{f(c)}$ admits a left  $L_K(E)$-module action (which we appropriately denote by $\cdot$) by setting 
$$ r\cdot u = \Phi(r) (u)$$
for all $r\in L_K(E)$ and all $u\in \widehat{U}_{f(c)}$.  Finally,  by again using Remark \ref{equalonBremark}, we have that  this $L_K(E)$-action on $\widehat{U}_{f(c)}$ extends the $L_K(E)$-action on $U_{f(c)}$, i.e.,  that $U_{f(c)}$ is an $L_K(E)$-submodule of $\widehat{U}_{f(c)}$.  
\end{proof}

{Referring to the Figure~\ref{Figure:ref}, the modules  $\widehat U_{w-1}$,  $\widehat U_{f(g_1g_2g_3)}$,  and $\widehat U_{f(\ell)}$ contain proper formal series (and thus properly contain $U_{w-1}$,  $ U_{f(g_1g_2g_3)}$,  and $U_{f(\ell)}$, respectively).  
On the other hand,  $\widehat U_{f(d_1d_2d_3d_4)}=U_{f(d_1d_2d_3d_4)}$.}

Modules of the form $\widehat U_{f(c)}$ will play a central role in our analysis, as these will be shown to be the injective envelopes of the Chen simple modules of the form $V^f_{[c^\infty]}$.
We first  check that  $\widehat U_{f(c)}$ is an essential extension of $V^f_{[c^\infty]}$.

\begin{proposition}\label{prop:KEw} Let $E$ be a finite graph {with disjoint cycles}. For any cycle $c$, the simple left $L_K(E)$-module $V^f_{[c^\infty]}$ is essential in $\widehat U_{f(c)}$.
\end{proposition}
\begin{proof}
We establish that $U_{f(c)}$ is essential in $\widehat U_{f(c)}$; the result will then follow directly from Proposition \ref{Ufcfinflengthuniserial}, which in particular implies that $V^f_{[c^\infty]}$ is essential in $U_{f(c)}$.

By definition, any element in $\widehat U_{f(c)}$ is of the form
\[\sum_{j\in F_1}\sum_{h=0}^{\cdeg f(x)-1}\left(\sum_{\gamma\in \mathcal P^E_c}k_{j,\gamma,h}\gamma\right)\overline x^h\alpha^{\overline x}_{c,j}\]
for a suitable (finite) subset $F_1$ of $\mathbb I_c$.
Assume this element is not equal to zero. Then there exist $ j_0\in F_1$, $\gamma_0\in\mathcal P^E_c$, and $0\leq h_0<\cdeg f(x)$ such that $k_{j_0,\gamma_0, h_0}\not=0$. 
By Remark \ref{anyclosedpath},  the paths in $\mathcal P^E_c$ do not end in any closed path.
So by Lemma~\ref{lemma:tecnico2} we have
\[\gamma_0^*\sum_{j \in F_1}\sum_{h=0}^{\cdeg f(x)-1}\sum_{\gamma\in\mathcal P^E_c}k_{j,\gamma}\gamma \overline x^h\alpha^{\overline x}_{c,j} \ = \ 
\sum_{j\in F_1} \sum_{h=0}^{\cdeg f(x)-1} k_{j,\gamma_0, h} \overline x^h\alpha^{\overline x}_{c,j}.\]
The latter is a nonzero element in $U_{f(c)}$ since by Lemma~\ref{lemma:base} the  set $\{ \overline x^h\alpha^{\overline x}_{c,j} \ | \ j\in F_1$, $0\leq h<\cdeg f(x)\}$ is $K$-linearly independent,  and $k_{j_0,\gamma_0, h_0}\not=0$. 
\end{proof}

\section{The main result}\label{mainsection}

The aim of this paper is to explicitly construct the injective envelope of all simple modules over any Leavitt path algebra $L_K(E)$ associated to a finite graph $E$ {with disjoint cycles}.  A first important step has already been established.   Using a number of powerful tools (e.g., pure injectivity, and the B\'{e}zout property of Leavitt path algebras), the three authors proved the following. 

\begin{theorem}\label{thm:amt19}\cite[Theorem~6.4]{AMT19}
Let $E$ be any finite graph, and let $c$ be a proper cycle in $E$.  Then the Pr\"ufer module $U_{E,c-1}$  is injective if and only if $c$ is a maximal cycle (i.e.,  there are no cycles in $E$ other than cyclic shifts of $c$ which connect to $s(c)$).
\end{theorem}
{Referring to the Figure~\ref{Figure:ref}, $U_{E,c-1}$ is injective if and only if $c=d_1d_2d_3d_4$.}

With Theorem \ref{thm:amt19} and the results  established in the previous sections in hand, we are now in position to state the main theorem of the article.  

\begin{theorem}\label{thm:main}
Let $K$ be any field, and let $E$ be any finite graph {with disjoint cycles}.
Let  $c$ be a  cycle in $E$, and let $f(x)\in K[x]$ be a basic irreducible polynomial in $K[x]$. Then
the injective envelope of $V^f_{[c^\infty]}$ is the left $L_K(E)$-module $\widehat U_{f(c)}$.\\
In particular, if $c=w$ is a sink, then the injective envelope of $V_{[w^\infty]}=L_K(E)w$ is the left $L_K(E)$-module
$\widehat U_{w-1}=K[[E]]w$.
\end{theorem}

Observe that for a proper cycle $c$,
 $\mathcal P^E_c$ is finite if and only if $c$ is a maximal cycle. In such a situation we get $U_{f(c)}= \widehat U_{f(c)}$.
Thus Theorem~\ref{thm:main} may be viewed as a significant extension of 
Theorem \ref{thm:amt19}
for graphs {with disjoint cycles}, as it provides the injective envelopes of simples associated to all sinks, and to all cycles (not only the maximal ones) together with appropriate irreducible polynomials of $K[x]$.  

\begin{remark}\label{injcogenremark}
As mentioned in the Introduction, by \cite[Theorem 1.1]{AR14},  when $E$ {has disjoint cycles} then the Chen simple modules represent {\rm all} the simple $L_K(E)$-modules.   So, once the proof of Theorem \ref{thm:main} is completed, we will get as a consequence an explicit description of an injective cogenerator for $L_K(E)$-Mod, namely, the direct product of the injective envelopes of all the simple modules. 
\end{remark}
\medskip
We will present the proof of Theorem \ref{thm:main}  below.   Here is an overview of how we will proceed.  In Proposition \ref{prop:KEw} we have already checked that $V_{[c^\infty]}^f$  is an essential submodule of $\widehat U_{f(c)}$. In Proposition~\ref{prop:inj1} we reduce the task of proving  the injectivity of $\widehat U_{f(c)}$ for any basic irreducible $f(x)\in K[x]$ to the specific polynomial $f(x)=x-1$. To do so, we will use that $\widehat U_{f(c)}$ is obtained from the left $L_{K'}(E)$-module $\widehat U_{\overline x^{-1}c-1}$ by applying the functor $\mathcal U$ which restricts the scalars from $K'$ to $K$.

%
%

Then,  to prove the injectivity of $\widehat U_{c-1}$, for any cycle $c$ in $E$, we reduce first to connected graphs,  then reduce to graphs  which contain no source vertices, and finally reduce to graphs   in which every source cycle is a loop. 

After these four reductions,   the core of the proof is based on an induction argument  on the number of cycles in the graph $E$. It is here that the main ideas come to light.

\begin{proposition}\label{prop:inj1}
Let $c$ be a proper cycle in $E$ and $f(x)\in K[x]$  a basic irreducible polynomial. Denote by $K'$ the field $K[x]/\langle f(x)\rangle$.  If $\widehat U_{c-1}$ is an injective left $L_{K'}(E)$-module, then $\widehat U_{f(c)}$ is an injective left $L_K(E)$-module.
\end{proposition}
\begin{proof}
Given $0\not=a\in K$, we first consider the case $f(x)=a^{-1}x-1$. In this situation, since ${\rm deg}(f(x))=1$ we get  $K=K'$. Then $\widehat U_{f(c)}$ is the  $L_K(E)$-module $\widehat U_{a^{-1}c-1}$, which is the twisted version of the $L_K(E)$-module $\widehat U_{c-1}$. Since twisting is an autoequivalence of $L_K(E)$-Mod,  $\widehat U_{a^{-1}c-1}$ is injective.

If $f(x)$ is a basic irreducible polynomial of degree $>1$, writing $f(x) = xg(x) - 1$ we get from the previous paragraph that $\widehat U_{g(\overline x)c-1}$ is an injective left $L_{K'}(E)$-module. Now consider the left $L_K(E)$-module $\widehat U_{f(c)}=\mathcal U\left(\widehat U_{g(\overline x)c-1}\right)$. By \cite[\href{https://stacks.math.columbia.edu/tag/015Y}{Lemma~12.29.1}]{stacks-project}, since $\mathcal U$ is right adjoint to 
\[-\otimes_{K}K':L_K(E)\text{-Mod}\to L_{K'}(E)\text{-Mod},
\]
and the latter is exact, $\mathcal U$ transforms injectives into injectives.
\end{proof}

The next three lemmas will permit us to reduce our study to connected graphs having no source vertices and in which the source cycles are loops. 

\begin{lemma}\label{lemma:connected}
Let $E$ be a finite graph and $E_i$, $i=1,2,..., m$, its connected components. The Leavitt path algebra associated to $E$ decomposes as the direct sum of two-sided ideals, each of which is the  Leavitt path algebra associated to a connected component of $E$:
\[L_K(E)=\bigoplus_{i=1}^m L_K(E_i).\]
For any cycle $c$ in $E_1$, we have
\[\widehat U_{E,c-1}=\widehat U_{E_1,c-1}.\]
Moreover, let $J$ be an ideal of $L_K(E)$ and $\rho_1$ be the sum of all vertices in $E_1$. We have $J=\rho_1J\oplus (1-\rho_1) J$ and $\varphi\left((1-\rho_1)J \right)=0$ for each morphism $\varphi:J\to  \widehat U_{E,c-1}=\widehat U_{E_1,c-1}$.
\end{lemma}
\begin{proof}
For the decomposition $L_K(E)=\bigoplus_{i=1}^m L_K(E_i)$ see \cite[Proposition~1.2.14]{AAS17}. 
Clearly $\mathcal P^E_c=\mathcal P^{E_1}_c$. Therefore we get $\widehat U_{E,c-1}=\widehat U_{E_1,c-1}$. Finally, for any
\[\varphi:J\to  \widehat U_{E,c-1}=\widehat U_{E_1,c-1}\]
we have 
\[\varphi\left((1-\rho_1)J\right)=(1-\rho_1)\varphi(J)\subseteq (1-\rho_1)\widehat U_{E_1,c-1}.\]
Since $(1-\rho_1)\mathcal P^{E_1}_c=0$ we get $\varphi\left((1-\rho_1)J\right)=0$.
\end{proof}

The reduction to graphs having no source vertices and in which the source cycles are loops will be achieved using relatively standard equivalence functors between categories of modules over Leavitt path algebras associated to  general graphs, and categories of modules over Leavitt path algebras associated to graphs having no source vertices  and all source cycles being loops.   
 Of course such equivalences will preserve various homological properties, including injectivity.  

However, we need more. The modules $ \widehat U_{E,c-1} $ are built in a specific way from the data corresponding to the graph $E$ and cycle $c$.   So if, for example, $F$ is a subgraph of $E$ (or $F$ is otherwise built from $E$)  which contains the cycle $c$, and if $T: L_K(F)\text{-Mod} \to L_K(E)\text{-Mod}$ is an equivalence of categories, and if $ \widehat U_{F,c-1}$ has been shown to be injective as an $L_K(F)$-module, then certainly $T(\widehat U_{F,c-1})$  is an injective $L_K(E)$-module.   But it is not at all immediate from purely categorical considerations that
\[T(\widehat U_{F,c-1}) \cong \widehat U_{E,c-1}.\]
Fortunately,  the  displayed isomorphism does indeed hold for each of the equivalence functors  that we will utilise.

 \begin{lemma}\label{lemma:nosources}
 Let $E$ be a finite graph {with disjoint cycles}, and $\overline u$ a source vertex in $E$ {(see Figure~\ref{Figure:ref2})}.   Let $H$ denote the hereditary subset $E^0 \setminus \{\overline{u} \}$ of $E^0$. 
  Setting $\varepsilon=\sum_{u\in E^0, u\not=\overline u}u$, there is  a Morita equivalence
 \[\xymatrix{L_K(E)\text{-Mod}\ar@<1ex>[rrr]^-{\varepsilon L_K(E)\otimes_{L_K(E)}-}&&&L_K(E_H)\text{-Mod}
 =\varepsilon L_K(E) \varepsilon\text{-Mod}
\ar@<1ex>[lll]^-{ L_K(E) \varepsilon\otimes_{\varepsilon L_K(E)\varepsilon}-}
}.\]
Moreover, for any cycle $c$ of $E$, necessarily $c$ is in $E_H$, and we have
\[\widehat U_{E,c-1}\ \cong\ L_K(E)\varepsilon\otimes_{\varepsilon L_K(E)\varepsilon} \widehat U_{E_H,c-1}.\]

 \end{lemma}
 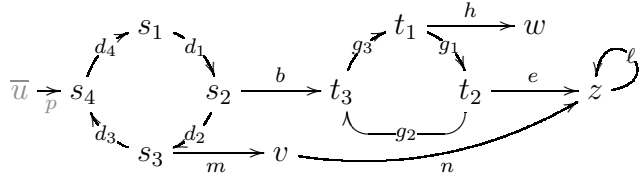
\begin{figure}
\begin{center}{
\[\xymatrix@-1.3pc{&&&s_1\ar@/^/[dr]|{d_1}&&&&t_1\ar@/^/[dr]|{g_1}\ar[rr]^h&&w\\\
&\textcolor{light-gray}{\overline u}\ar@[light-gray][r]_{\textcolor{light-gray}{ p}}&s_4\ar@/^/[ur]|{d_4}&&s_2\ar@/^/[dl]|{d_2}\ar[rr]^b&&t_3\ar@/^/[ur]|{g_3}&&t_2\ar@<-2pt> `d[l] `[ll]|{g_2} [ll]\ar[rr]^e&&z\ar@(r,u)|\ell\\
&&&s_3\ar@/^/[ul]|{d_3}\ar[rr]_m&&v\ar@/_1pc/ [rrrrru]_n}\]
}
\caption{Eliminating source vertices}
\label{Figure:ref2}
\end{center}
\end{figure}
 \begin{proof}
 That the indicated functors give a Morita equivalence has been proved in \cite[Lemma~4.3]{AR14}.
Let $c$ be a cycle in $E$. 
So
\[\mathcal P^E_c=\varepsilon\mathcal P^E_c\cup \overline u \mathcal P^E_c=
\mathcal P^{E_H}_c\cup
\overline u \mathcal P^E_c.\]
Any element $\gamma$ in $\overline u \mathcal P^E_c$ has the form
$\gamma=\gamma_1\varepsilon\gamma_2\cdots \gamma_\ell$, where $s(\gamma_1) = \overline{u}$ and $r(\gamma_1) \in H$.   
The map 
\[{K\mathcal P^E_c \alpha_{c,j} \ \to \ 
 L_K(E)\varepsilon\otimes_{\varepsilon L_K(E)\varepsilon}K\varepsilon \mathcal P^E_c\alpha_{c,j} = 
L_K(E)\varepsilon\otimes_{\varepsilon L_K(E)\varepsilon} K\mathcal P^{E_H}_c \alpha_{c,j}}
\] defined by $\gamma  \alpha_{c,j}\mapsto \varepsilon \otimes\gamma  \alpha_{c,j}$ for each $\gamma\in{\varepsilon\mathcal P^E_c}$, and $\gamma \alpha_{c,j}\mapsto \gamma_1\otimes \gamma_2\cdots  \gamma_\ell\alpha_{c,j}$ for each $\gamma=\gamma_1\varepsilon\gamma_2\cdots \gamma_\ell\in  \overline u  \mathcal P^E_c$, induces an isomorphism
\[\xymatrix{{}\widehat U_{E,c-1}\ar[rr]&& L_K(E)\varepsilon\otimes_{\varepsilon L_K(E)\varepsilon}\widehat U_{E_H,c-1}\\
{}\bigoplus_{j\in\mathbb I_c}K[[\mathcal P^E_c]]\alpha_{c,j}\ar@{=}[u]&&
L_K(E)\varepsilon\otimes_{\varepsilon L_K(E)\varepsilon} \left( \bigoplus_{j\in\mathbb I_c}K[[\mathcal P^{E_H}_c]]\alpha_{c,j} \right) \ar@{=}[u]}\]
{(where the direct sums are as $K$-vector spaces),} thus establishing the result.
 \end{proof}
 
 Assume now that 
$E$ contains  a source cycle $d=d_1\cdots d_r$ which is not a loop. 
By \cite[Lemma~4.4]{AR14} a finite graph $F_{E,d}$ can be constructed from $E$ in which the cycle $d$ (and all its vertices) is replaced by a loop in such a way that $L_K(E)$ and $L_K(F_{E,d})$ are Morita equivalent. More precisely,  let $d^0$ denote the set of vertices $\{s(d_1), \dots , s(d_r)\}$, and define  $F_{E,d}^0=\{\overline v\}\cup (E^0\setminus d^0)$ where $\overline v$ is a new vertex. Then to define $F_{E,d}^1$, we first set $s_{F_{E,d}}^{-1}(u)=s_E^{-1}(u)$ for each $u\in E^0\setminus d^0$; for each edge $f$ with $s(f)\in d^0$ and $r(f)\in E^0\setminus d^0$, define an edge $\varphi(f)$ with $s_{F_{E,d}}(\varphi(f))=\overline v$ and $r_{F_{E,d}}(\varphi(f))=r_E(f)$. Finally, we define a loop $d'$ at $\overline v$ so that $s_{F_{E,d}}(d')=\overline v=r_{F_{E,d}}(d')$. 
Observe that there are no edges connecting any pair of  vertices in $d^0$ other than the edges  $d_1$, ..., $d_r$.

 \begin{lemma}\label{lemma:sourcecycle}
 Let $d=d_1\cdots d_r$ be a source cycle in the graph $E$ with $r\geq 2$.  Let $F_{E,d}$ be the graph described in the previous paragraph, in which $d$ is replaced by the loop $d'$.  Then  there are inverse equivalence functors
  \[\xymatrix{L_K(E)\text{-Mod}\ar@<1ex>[rrr]^{G_1}&&&L_K(F_{E,d})\text{-Mod}\ar@<1ex>[lll]^{G_2}}.\]
  Moreover,  $
\widehat U_{E,c-1}\cong G_2\left( \widehat U_{F_{E,d},c-1}\right)$ for any cycle $c$ in $E$.  
%
%
%
 \end{lemma}

 \begin{proof}
We describe   how the functors $G_1$ and $G_2$ work.  (The details are given in \cite[Lemma 4.4]{AR14}.)  
Consider the map $\theta: L_K(F_{E,d})\to L_K(E)$
defined on vertices by: 

\medskip

$\theta(u)=u \ \text{ for }u\in F_{E,d}^0\setminus \{\overline v\};$ 

\smallskip

$ \theta(\overline v)=s_E(d);$

\medskip

\noindent
and defined on edges by:

\medskip

$\theta(e)=e\text{ for }e\text{ having }s_{F_{E,d}}(e)\in F_{E,d}^0\setminus\{\overline v\}$;

\smallskip

$\theta(\varphi(f))=d_1\cdots d_{i-1}f\text{ for }f \text{ having  }s_E(f)=s(d_i)\text{ and }r_E(f)\in E^0\setminus d^0;$

\smallskip


$\theta({d'})=d=d_1\cdots d_r.$  

\medskip

\noindent
By \cite[Lemma~4.4]{AR14} the map $\theta$ extends to a well-defined $K$-algebra isomorphism
$$\theta:  L_K(F_{E,d}) \to \omega L_K(E)\omega,$$ 
where $\omega:=s(d_1)+\sum_{u\in E^0\setminus d^0}u$.  We denote by $\hat\theta$ the {induced} natural isomorphism of categories
$$  \hat\theta:  L_K(F_{E,d})\text{-Mod} \to \omega L_K(E)\omega\text{-Mod}.$$
 The ring $\omega L_K(E)\omega$ is easily seen to be a full corner of $L_K(E)$ and so  is Morita equivalent to $L_K(E)$.
The functors $G_1$ and $G_2$ which realise the Morita equivalence between $L_K(E)$-Mod and $L_K(F_{E,d})$-Mod are given by the following compositions
\[\xymatrix{L_K(E)\text{-Mod}\ar@<1ex>[rrr]^{\omega L_K(E)\otimes_{L_K(E)}-}&&&\omega L_K(E) \omega\text{-Mod}\ar@<1ex>[r]^-{\hat\theta^{-1}}\ar@<1ex>[lll]^{ L_K(E) \omega\otimes_{\omega L_K(E)\omega}-}&L_K(F_{E,d})\text{-Mod}\ar@<1ex>[l]^-{\hat\theta}}.\]
%
Now let  $c$ be a  cycle in $E$ (possibly equal to $d$).   Then
\[\mathcal P^E_c=\omega \mathcal P^E_c\cup (1-\omega)\mathcal P^E_c.\]
The map
\[U_{E,c-1}=\bigoplus_{j\in\mathbb I_c}K\mathcal P^E_c \alpha_{c,j} \ \ \to \ \  L_K(E)\omega\otimes_{\omega L_K(E)\omega}\left(\bigoplus_{j\in\mathbb I_c}K \omega \mathcal P^E_c \alpha_{c,j}\right)\]
defined by
$$\rho \alpha_{c,j} \ \mapsto \  \omega\otimes \rho \alpha_{c,j}  \ \   \mbox{ for each } \rho\in \omega \mathcal P^E_c   $$
\noindent
 and by 
\[s(d_i)p\alpha_{c,j} \ \mapsto  \ d^*_{i-1}\cdots d_1^*\otimes {d_1\cdots d_{i-1}s(d_i)p}\alpha_{c,j}  \ \   \mbox{ for each } s(d_i)p\in (1-\omega)\mathcal P^E_c    \]
 extends to an isomorphism 
 $$\widehat{U}_{E,c-1}=\bigoplus_{j\in\mathbb I_c}K[[\mathcal P^E_c]]\alpha_{c,j}  \ \ \to  \ \ 
 L_K(E)\omega\otimes_{\omega L_K(E)\omega} \left(\bigoplus_{j\in\mathbb I_c}K[[\omega \mathcal P^E_c]]\alpha_{c,j}\right).$$ 
Since $\theta(\mathcal P^{F_{E,d}}_c)=\omega \mathcal P^E_c$, we have
$\hat\theta(\widehat U_{F_{E,d},c-1})\cong\bigoplus_{j\in\mathbb I_c}K[[\omega \mathcal P^E_c]]\alpha_{c,j}$, and hence
\[G_2\left(\widehat U_{F_{E,d},c-1}\right)=L_K(E)\omega\otimes_{\omega L_K(E)\omega}\hat\theta(\widehat U_{F_{E,d},c-1})\cong
\widehat U_{E,c-1},
\]
as desired.
 \end{proof}
 
 In the following example we  clarify in a concrete situation how the isomorphisms described in Lemma \ref{lemma:sourcecycle} work.
 
 \begin{example}
 Let $E'$ be the graph
\[\xymatrix@-1.3pc{&&\textcolor{light-gray}{s_1}\ar@/^/[dr]|{\textcolor{light-gray}{d_1}}&&&&t_1\ar@/^/[dr]|{g_1}\ar[rr]^h&&w\\\
&\textcolor{light-gray}{s_4}\ar@/^/[ur]|{\textcolor{light-gray}{d_4}}&&\textcolor{light-gray}{s_2}\ar@/^/[dl]|{\textcolor{light-gray}{d_2}}\ar[rr]^{\textcolor{light-gray}{b}}&&t_3\ar@/^/[ur]|{g_3}&&t_2\ar@<-2pt> `d[l] `[ll]|{g_2} [ll]\ar[rr]^e&&z\ar@(r,u)|\ell\\
&&\textcolor{light-gray}{s_3}\ar@/^/[ul]|{\textcolor{light-gray}{d_3}}\ar[rr]_{\textcolor{light-gray}{m}}&&v\ar@/_1pc/ [rrrrru]_n}\]
%
and
\[\xymatrix@-1.3pc{&&&&t_1\ar@/^/[dr]|{g_1}\ar[rr]^h&&w\\\
F_{E',d}=&\bm{\overline v}\ar@(l,u)|{\bm{d'}}\ar[rr]^{\bm{\varphi(b)}}\ar[rd]_{\bm{\varphi(m)}}&&t_3\ar@/^/[ur]|{g_3}&&t_2\ar@<-2pt> `d[l] `[ll]|{g_2} [ll]\ar[rr]^e&&z\ar@(r,u)|\ell\\
&&v\ar@/_1pc/ [rrrrru]_n}\]
the graph in which the cycle $d_1d_2d_3d_4$ is substituted by the loop $d'$. Denote by $\omega$ the idempotent $s_1+\sum_{u\in {E'}^0\setminus d^0}u$ of $L_K(E')$, and by $c$ either one of the the cycles $g_1g_2g_3$  and  $\ell$, or the sink $w$. 
In the isomorphism 
  \[\widehat U_{E',c-1}\cong L_K(E')\omega\otimes_{\omega L_K(E')\omega}\hat\theta(\widehat U_{F_{E',d},c-1}),\]
{for any $\overline j\in\mathbb I_c$} the element $$\sum_{i=0}^\infty d_4d^id_1b\gamma_i{\alpha_{c,\overline j}}$$ of $\widehat U_{E',c-1}$ corresponds to the element 
 \[d_3^*d_2^*d_1^*\otimes \sum_{i=0}^\infty \theta\left((d')^{i+1}\varphi(b)\gamma_i \right){\alpha_{c,\overline j}}=d_3^*d_2^*d_1^*\otimes \sum_{i=0}^\infty d^{i+1}d_1b\gamma_i{\alpha_{c,\overline j}}\]
 of  $L_K(E')\omega\otimes_{\omega L_K(E')\omega} \left(\bigoplus_{j\geq 1}K[[\omega \mathcal P^{E'}_c]]\alpha_{c,j}\right)$.  
%
%
%
%
 \end{example}
 \medskip

  \bigskip

We are now ready to give the 
\begin{proof}[Proof of Theorem~\ref{thm:main}]
{Recall that we consider} any sink $w\in E^0$ as a cycle of length 0. 

The essentiality of $V^f_{[c^\infty]}$ in $\widehat U_{f(c)}$ has been shown in Proposition~\ref{prop:KEw}.  So we need only demonstrate that $\widehat U_{f(c)}$ is injective  for each cycle $c$ and each basic irreducible polynomial $f(x)\in K[x]$.

Given any cycle $c$, by Lemma~\ref{lemma:connected} we can establish the injectivity of $\widehat U_{f(c)}$ in the Leavitt path algebra associated to the connected component containing $c$.
By Lemmas~\ref{lemma:nosources} and \ref{lemma:sourcecycle},  we can assume that the connected component containing $c$ has no source vertices,  and that  all its source cycles are source loops. Moreover by Proposition~\ref{prop:inj1},  it suffices to show the result in the case $f(x)=x-1$, and hence we may assume  $\widehat U_{f(c)}= \widehat U_{c-1}$.
\begin{figure}
\begin{center}
\[\xymatrix@-1.3pc{&&&&t_1\ar@/^/[dr]|{g_1}\ar[rr]^h&&w\\\
&{t}\ar@(l,u)|{{\tau}}\ar[rr]^{{\varepsilon_1}}\ar[rd]_{{\varepsilon_2}}&&t_3\ar@/^/[ur]|{g_3}&&t_2\ar@<-2pt> `d[l] `[ll]|{g_2} [ll]\ar[rr]^e&&z\ar@(r,u)|\ell\\
&&v\ar@/_1pc/ [rrrrru]_n}\]
\caption{Graph with disjoint cycles, without source vertices, whose source cycle is a loop}
\label{Figure:ref3}
\end{center}
\end{figure}
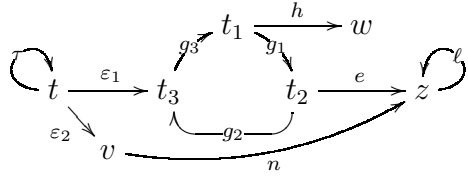
\smallskip

We proceed by induction on the number $N$ of cycles in  the connected component of $E$ containing the cycle $c$. We continue to call this connected component $E$.

If $N=0$, then by the reduction assumptions the graph $E$ has one vertex and no edges. So $L_K(E)\cong K$, and therefore {\it all}  left $L_K(E)$-modules are injective. 

{If $N=1$, either $E$ consists of one vertex $u$ and one loop at $u$ or, invoking Remark \ref{rem:rigidity} {and the assumption that $N=1$} the graph $E$ contains  
 {a} source loop $\tau$ and  at least two vertices. In the first case by Theorem~\ref{thm:amt19} we conclude that   $U_{c-1} =  \widehat U_{c-1}$ is injective. In the second case, since
$\tau$ is a maximal cycle, by Theorem \ref{thm:amt19} the Pr\"ufer module $U_{E,\tau-1}=\widehat U_{E,\tau-1}$ is injective.}

{Now let $N>1$. {If } $\tau$ is a maximal cycle, then  again by Theorem \ref{thm:amt19} the Pr\"ufer module $U_{E,\tau-1}=\widehat U_{E,\tau-1}$ is injective. So we have to prove that $\widehat U_{E,c-1}$ is injective for any non-maximal cycle $c$. Denote $s(\tau)$ by $t$.
Denote by $I$ the two sided ideal of $L_K(E)$ generated by the hereditary set of vertices $H = E^0\setminus\{t\}$, and by $\rho\in I$ the idempotent $\sum_{u\in H}u {=1_{L_K(E)}-t}$. Moreover, let $s^{-1}(t)=\{\varepsilon_1,...,\varepsilon_n\}$ (see, e.g., Figure~\ref{Figure:ref3}).
}
%


By Proposition~\ref{lemma:idealeprincipale}, 
$I$ is isomorphic to a Leavitt path algebra, and there is an equivalence  between the category $I$-Mod of (unitary) $I$-modules and the category $\rho L_K(E)\rho$-Mod$=\rho I\rho$-Mod,  induced by the functors 
\[\rho I\otimes_I-:I\text{-Mod}\to \rho  L_K(E)\rho\text{-Mod}\quad \ \ \text{and}\]
\[\quad I\rho\otimes_{\rho L_K(E)\rho}-:\rho L_K(E)\rho\text{-Mod}\to I\text{-Mod}.\]
Moreover, the $K$-algebra $\rho L_K(E)\rho$ is isomorphic to the Leavitt path $K$-algebra $L_K(E_H)$. 

Now let $c$ be any cycle in $E$, $c \neq \tau$.   Then $c$ lives also in $E_H$.
Since each of the connected components of the graph $E_H$ (there might be more than one such component) has fewer cycles  than does $E$, we can apply the inductive hypothesis to conclude that  $\widehat U_{E_H,c-1}$ is an injective $L_K(E_H)$-module.

We prove that $\widehat U_{E,c-1}$ is injective
in $L_K(E)$-Mod. To do so we use  the version of Baer's Criterion  presented in Lemma~\ref{injectivebystep}, with respect to the ideal $I$.

So first, let $J$ be a left ideal of $L_K(E)$ contained in $I$. Since $I$ has local units, we have $IJ=J$ and hence $J$ belongs to $I$-Mod. Consider an arbitrary $L_K(E)$-homomorphism
\[\varphi:J\to \widehat U_{E,c-1}.\] 
Since $IJ=J$, $\varphi$ factors through 
\[\overline{\varphi}:J\to I \widehat U_{E,c-1}.\]
We show  that it is possible to extend $\varphi$ to $I$. Clearly $I\widehat U_{E,c-1}$ is also a left $I$-module and $\overline{\varphi}$ is an $I$-linear map.
Applying the functor $\rho I\otimes_I-$ we get the $\rho L_K(E)\rho$-homomorphism
\[\rho I\otimes_I \overline{\varphi}: \ \ \rho I\otimes_I J \ \to \ \rho I\otimes_I I\widehat U_{E,c-1}.\]

{{\bf Claim}:   $\rho I\otimes_I I\widehat U_{E,c-1} \cong \widehat U_{E_H,c-1}$ as left $\rho L_K(E)\rho \cong L_K(E_H)$-modules.   } 
{{\bf Proof}. Any element in $\rho I\otimes_I I\widehat U_{E,c-1}$ is of the form
$\rho \iota_1\otimes \iota_2z$ where $\iota_1,\iota_2\in I$ and $z\in\widehat U_{E,c-1}$. Denoting by $\iota$ the product $\iota_1\iota_2$ we have
\[\rho \iota_1\otimes \iota_2z=\rho^2 \iota_1\otimes \iota_2z=\rho \otimes \rho \iota_1\iota_2 z=\rho \otimes \rho \iota z,\]
and $\rho\iota z$ has the form
\[\sum_{j\in F_1}\sum_{\gamma\in \mathcal P^{E_H}_c} k_{j,\gamma}\rho\iota\gamma \alpha_{c,j} +\sum_{j\in F_2}\sum_{\gamma\in \mathcal P^{E_H}_c}\sum_{y=1}^n \sum_{x=0}^\infty  k_{j,\gamma,x,y }\rho\iota\tau^x\varepsilon_y\gamma \alpha_{c,j}.\]
}
{We establish now the key point, which is that  if  $ \iota \in I$ and $ z \in I\widehat U_{E,c-1}$ then  $\rho \iota z \in \widehat U_{E_H,c-1}$.    
Note that for each $\mu\in \mathcal P^E_c$ we have
\[\rho\mu=\begin{cases}
    \mu  & \text{if $\mu$ is a path in  $\mathcal P^{E_H}_c$}, \\
     0 & \text{if $s(\mu)=t$}.
\end{cases}
\]
Separating the paths in $\mathcal P^{E}_c$ starting from the vertex $t$ and those starting from a vertex in $H$, any element $z$ of $\widehat U_{E,c-1}$ has the form (see, e.g., Figure~\ref{Figure:ref3})
\[z=\sum_{j\in F_1}\sum_{\gamma\in \mathcal P^{E_H}_c} k_{j,\gamma}\gamma \alpha_{c,j} \  + \  \sum_{j\in F_2}\sum_{\gamma\in \mathcal P^{E_H}_c}\sum_{y=1}^n \sum_{x=0}^\infty  k_{j,\gamma,x,y }\tau^x\varepsilon_y\gamma \alpha_{c,j}\]
for suitable finite subsets $F_1$ and $F_2$ of $\mathbb I_c$ and $k_{j,\gamma}$, $k_{j,\gamma,x,y }\in K$. Then 
\[\rho \iota z = \sum_{j\in F_1}\sum_{\gamma\in \mathcal P^{E_H}_c} k_{j,\gamma}\rho\iota\gamma \alpha_{c,j}  \ + \ \sum_{j\in F_2}\sum_{\gamma\in \mathcal P^{E_H}_c}\sum_{y=1}^n \sum_{x=0}^\infty  k_{j,\gamma,x,y }\rho\iota\tau^x\varepsilon_y\gamma \alpha_{c,j}.\]
The left hand double summand   is an element of $\widehat U_{E_H,c-1}$. In the  right hand quadruple summand  consider each $\rho\iota\tau^x \varepsilon_y \gamma \alpha_{c,j}$. The term $\rho\iota$ is a linear combination of monomials in $L_K(E)$ starting from vertices in $H$,  and $\gamma$ is an element of $\mathcal P^{E_H}_c$. 
 In the linear combination of monomials $\rho\iota$, let $\rho \iota'_{x,y}\varepsilon_y^*(\tau^*)^x$ be the sub-linear combination of those monomials  ending with $\varepsilon_y^*(\tau^*)^x$. Since $\iota$ is fixed, we have $\iota'_{x,y}\not=0$ for a finite number of $(x,y)\in \mathbb N\times\{1,...,n\}$, and, in particular,  $\iota'_{x,y}=0$ for $x$ greater than a suitable $\ell\in\mathbb N$. Therefore
\[\sum_{j\in F_2}\sum_{\gamma\in \mathcal P^{E_H}_c}\sum_{y=1}^n \sum_{x=0}^\infty  k_{j,\gamma,x,y }\rho\iota\tau^x\varepsilon_y\gamma \alpha_{c,j}=
\sum_{j\in F_2}\sum_{\gamma\in \mathcal P^{E_H}_c}\sum_{y=1}^n \sum_{x=0}^\ell  k_{j,\gamma,x,y }\rho\iota'_{x,y}\rho\gamma  \alpha_{c,j},\]
which is an element of $\widehat U_{E_H,c-1}$.    Thus  $\rho \iota z \in \widehat U_{E_H,c-1}$ as desired. }

{    Now consider the map $(\rho \iota_1, \iota_2 z) \mapsto  \rho \iota_1 \iota_2 z$ where $\iota_1, \iota_2 \in I$ and $z \in \widehat U_{E,c-1}$.   Clearly this is an $I$-balanced  linear map, which with the previously established key point yields a left $\rho L_K(E) \rho$-module map  $\psi:  \rho I\otimes_I I\widehat U_{E,c-1} \to  \widehat U_{E_H,c-1}$.      
But by the induction hypothesis, since $E_H$ has fewer cycles than does $E$, $\widehat U_{E_H,c-1}$ is an 
 an injective module over the ring $\rho I\rho=\rho L_K(E)\rho \cong L_K(E_H)$.
}

Since $J\leq I$ and $IJ=J$, applying the exact functor $\rho I\otimes_ I-$ 
    we get 
the solid part of the following commutative diagram of left $\rho L_K(E)\rho$-modules:
{\[\xymatrix{0\ar[r]&\rho I\otimes_I J\ar[rr]\ar[d]_{\rho I\otimes \overline{\varphi}}&&\rho I\otimes_I I\ar@{-->}^\psi[lld]\\
&\rho I\otimes_I I\widehat U_{E,c-1} \cong\widehat U_{E_H,c-1}}.\]}
The existence of the dashed arrow $\psi$  follows by the injectivity of $\widehat U_{E_H,c-1}$ in $\rho L_K(E)\rho$-Mod established above.   Applying the equivalence functor $I\rho\otimes_{\rho L_K(E)\rho}-$, 
we get the following commutative diagram of $I$-modules:
{\[\xymatrix@-1pc{I\rho\otimes_{\rho L_K(E)\rho}(\rho I\otimes_I J)\cong J\ar[r]\ar[dd]_{I\rho\otimes\left(\rho I\otimes \overline{\varphi}\right)\cong \overline{\varphi}}& I\rho\otimes_{\rho L_K(E)\rho}(\rho I\otimes_I I)\cong I\ar[ldd]^{\overline{\psi}}
\\
\\
I\rho\otimes_{\rho L_K(E)\rho} (\rho I\otimes_I I\widehat U_{E,c-1})\cong I\widehat U_{E,c-1}
}\]}
Moreover, the $I$-linear map $\overline{\psi}$ is in fact $L_K(E)$-linear. Indeed, since $I$ has local units, for each $\iota\in I$ there exists $\zeta\in I$ such that $\zeta \iota=\iota$. Then for each $\lambda\in L_K(E)$, since $\lambda\zeta$ belongs to $I$, we have
\[\overline\psi(\lambda \iota)=\overline\psi((\lambda \zeta)\iota)=  (\lambda\zeta)\overline\psi(\iota)=\lambda \cdot \overline\psi(\zeta\iota)=\lambda \cdot \overline\psi(\iota).\]
Composing $\overline{\psi}$ with the inclusion of $I\widehat U_{E,c-1}$ inside $\widehat U_{E,c-1}$ we get the desired extension of $\varphi$.   This completes the first step of the application of Baer's Criterion Lemma~\ref{injectivebystep}.  

\medskip

To complete  the second step required to apply  Lemma~\ref{injectivebystep}, we have to prove that any morphism of left $L_K(E)$-modules from a left ideal $J'$ containing $I$
extends to $L_K(E)$.   As before, we denote    $s^{-1}(t)$ by $\{\tau, \varepsilon_1, ..., \varepsilon_n\}.$

\medskip

We first assume $J'=I$. Let
\[\chi: I\to \widehat U_{E,c-1}\]
be a morphism of $L_K(E)$-modules. 
We have to extend $\chi$ to a morphism
$\hat\chi:L_K(E)=I+L_K(E)t\to \widehat U_{E,c-1}$.
The elements of the intersection $I\cap L_K(E)t$ have the following form:
\[\sum_{i,j} \ell_{i,j}(\varepsilon_i)^*(\tau^*)^{j}\quad\text{with}\quad \ell_{i,j}\in L_K(E).\]
In particular the restriction of $\chi$ to $I\cap L_K(E)t$ is determined by
$\chi\big((\varepsilon_i)^*(\tau^*)^{j}\big)$, $i\in\{1,...,n\}$, $j\in\mathbb N$. Observe that
\[z_{i,j}:=\chi\big((\varepsilon_i)^*(\tau^*)^{j}\big)=\chi\big(\rho(\varepsilon_i)^*(\tau^*)^{j}\big)=\rho\chi\big((\varepsilon_i)^*(\tau^*)^{j}\big)\in \rho \widehat U_{E,c-1}
.\]

Setting
\[
t \ \mapsto \ \sum_{i=0}^\infty\sum_{j=1}^n t\tau^i\varepsilon_jz_{i,j}=\sum_{i=0}^\infty\sum_{j=1}^n \tau^i\varepsilon_jz_{i,j} \ ,\]
we define a map $\overline {\chi}:L_K(E)t\to \widehat U_{E,c-1}$ whose restriction to $I\cap L_K(E)t$ coincides with $\chi$.  Specifically, 
\begin{eqnarray*}
\overline {\chi}\big((\varepsilon_{i_0})^*(\tau^*)^{j_0}\big) & \  = \  &
\overline {\chi}\big((\varepsilon_{i_0})^*(\tau^*)^{j_0}t\big) \ = \ 
\big((\varepsilon_{i_0})^*(\tau^*)^{j_0}\big)\overline {\chi}(t) \\
& \ = \ & \big((\varepsilon_{i_0})^*(\tau^*)^{j_0}\big)\sum_{i=0}^\infty\sum_{j=1}^n \tau^i\varepsilon_jz_{i,j} \ = \ z_{i_0,j_0}.
\end{eqnarray*}
Thus we can define
\[\hat \chi: I+L_K(E)t=L_K(E)\to \widehat U_{E,c-1}\]
by setting
$\hat \chi(\iota+\lambda t):=\chi(\iota)+\overline {\chi}(\lambda t)$.  This completes the verification of the second step,  in the specific  case $J' = I$.


\medskip

Assume now $J'$ is a left ideal of $L_K(E)$ properly containing $I$. By Proposition \ref{lemma:idealeprincipale}, $J'$ is equal to $L_K(E)p(\tau)$ for a suitable polynomial $p(x)\in K[x]$ with $p(0)=1$.
Consider a map $\theta:J'=L_K(E)p(\tau)\to \widehat U_{E,c-1}$.

Denote by $P(x)$ the formal series in $K[[x]]$ such that $p(x)P(x)=1$, and write $P(x) = \sum_{i=0}^\infty h_ix^i$.  {As observed previously, any element $z$ of $\widehat U_{E,c-1}$ has the form
\[z=\sum_{j\in F_1}\sum_{\gamma\in \rho\mathcal P^{E}_c} k_{j,\gamma}\gamma \alpha_{c,j} +\sum_{j\in F_2}\sum_{\gamma\in \rho\mathcal P^{E}_c}\sum_{y=1}^n \sum_{x=0}^\infty  k_{j,\gamma,x,y }\tau^x\varepsilon_y\gamma \alpha_{c,j}.\]
Therefore
\[P(\tau)z=\big(\sum_{i=0}^\infty h_i\tau^i\big)z=\]
\[=\sum_{i=0}^\infty\sum_{j\in F_1}\sum_{\gamma\in \rho\mathcal P^{E}_c}k_i k_{j,\gamma}\tau^i\gamma \alpha_{c,j} +\sum_{i=0}^\infty\sum_{j\in F_2}\sum_{\gamma\in \rho\mathcal P^{E}_c}\sum_{y=1}^n \sum_{x=0}^\infty  k_ik_{j,\gamma,x,y }\tau^{x+i}\varepsilon_y\gamma \alpha_{c,j}\]
is again an element of $\widehat U_{E,c-1}$.
}

Setting $\pi(1)=P(\tau)\theta(p(\tau))$ we get a map $\pi:L_K(E)\to \widehat U_{E,c-1}$ which extends $\theta$.   This completes the verification of the second step, in the case $I \subsetneqq J'$.  

\medskip

Therefore, by Lemma~\ref{injectivebystep}, we conclude that $\widehat U_{E,c-1}$ is an injective $L_K(E)$-module,  thereby completing  the proof of Theorem \ref{thm:main}.   
\end{proof}


%
%
%
%


\begin{thebibliography}{10}
\providecommand{\url}[1]{{#1}}
\providecommand{\urlprefix}{URL }
\expandafter\ifx\csname urlstyle\endcsname\relax
  \providecommand{\doi}[1]{DOI~\discretionary{}{}{}#1}\else
  \providecommand{\doi}{DOI~\discretionary{}{}{}\begingroup
 \urlstyle{rm}\Url}\fi


 
 
 

 

 
 

 

  \bibitem{Ab83} G. Abrams, \emph{Morita equivalence for rings with local units}, Comm. Algebra {\bf 11} (1983), pp. 801-837. DOI: 10.1080/00927878308822881
  
\bibitem{AAS17}   G. Abrams, P. Ara, M. Siles Molina.  Leavitt path algebras.  Lecture Notes in Mathematics vol. 2191.   Springer Verlag, London, 2017.  ISBN-13:  978-1-4471-7344-1.  DOI:  10.1007/978-1-4471-7344-1

\bibitem{AMT19} G. Abrams, F. Mantese, A. Tonolo, \emph{Pr\"{u}fer modules over Leavitt path algebras}, J. Algebra App.  {\bf 18} (2019), 1950154 (28 pp.).   
DOI: 10.1142/S0219498819501548


\bibitem{AMT21} G. Abrams, F. Mantese, A. Tonolo, \emph{Injective modules over the Jacobson algebra $K\langle X,Y\mid XY=1\rangle$},  Canadian Bulletin of Mathematics {\bf 64}(2) (2021), pp 323--339. DOI: 10.4153/S0008439520000478


\bibitem{AMT15} G. Abrams, F. Mantese, A. Tonolo, \emph{Extensions of simple modules over Leavitt path algebras}, Journal of Algebra {\bf 431}  (2015), pp. 78 -- 106.  
DOI: 10.1016/j.jalgebra.2015.01.034 


\bibitem{AAJZ} A. Alahmadi, H. Alsulami, S.K. Jain, E. Zelmanov, \emph{  Leavitt path algebras of finite Gelfand-Kirillov dimension}, J. Algebra Appl. {\bf 11}(6) (2012), 1250225 (6 pp.) 
DOI: 10.1142/S0219498812502258



\bibitem{AR14} P. Ara, K. M. Rangaswamy, \emph{Finitely presented simple modules over Leavitt path algebras},  Journal of Algebra {\bf 417}  (2014), pp. 333 -- 352.  DOI: 10.1016/j.jalgebra.2014.06.032

\bibitem{B}  G. Bergman, \emph{Coproducts and some universal ring constructions}, Trans. A.M.S. {\bf 200} (1974), pp. 33--88.  DOI: 10.1090/s0002-9947-1974-0357503-7 

\bibitem{Ch12} X. W. Chen, \emph{Irreducible representations of Leavitt path algebras}, Forum Math. 20 (2012).
DOI: 10.1515/forum-2012-0020


 \bibitem{G} L. Gerritzen, \emph{Modules over the algebra of the noncommutative equation $yx = 1$}, Arch.
Math. {\bf 75} (2000), pp. 98--112.   DOI: 10.1007/pl00000437


\bibitem{HSV22} R. Hazrat, A. N. Sebandal, J. P. Vilela, \emph{Graphs with disjoint cycles classification via the talented
              monoid}, J. Algebra {\bf 593} (2022), pp. 319--340. DOI: 10.1016/j.jalgebra.2021.11.022
              
              
 \bibitem{I} M. Iovanov and A. Sistko, \emph{On the Toeplitz-Jacobson algebra and direct finiteness}, pp. 113--124 in \emph{Groups, Rings, Group Rings, and Hopf Algebras}, Contemporary Math {\bf 668}, Amer. Math. Soc., Providence, RI, 2017.    DOI: 10.1090/conm/688/13830
 
  \bibitem{J}   N. Jacobson, \emph{Some Remarks on One-Sided Inverses}, Proc. A.M.S. {\bf 1} (1950), pp. 352--355.  DOI: 10.1007/978-1-4612-3694-8\_6





\bibitem{stacks-project} {The {Stacks project authors}},
 \emph {The Stacks project},
  {\url{https://stacks.math.columbia.edu}},
 (2022).
\end{thebibliography}
\end{document}